\pdfoutput=1
\documentclass[opre,nonblindrev]{informs3-generic}

\DoubleSpacedXI 

\usepackage{endnotes}
\let\footnote=\endnote
\let\enotesize=\normalsize
\def\notesname{Endnotes}
\def\makeenmark{$^{\theenmark}$}
\def\enoteformat{\rightskip0pt\leftskip0pt\parindent=1.75em
  \leavevmode\llap{\theenmark.\enskip}}

\usepackage{natbib}
 \bibpunct[, ]{(}{)}{,}{a}{}{,}
 \def\bibfont{\small}

\usepackage{amsmath,amsfonts,amssymb}
\usepackage{algorithm}
\usepackage{algpseudocode}
\usepackage{xcolor}
\usepackage{bbm}

\usepackage[T1]{fontenc}

\newtheorem{heuristic}{\small\sc Method}
\newenvironment{varalgorithm}[1]
{\algorithm\renewcommand{\thealgorithm}{#1}}
{\endalgorithm}

\DeclareMathOperator*{\argmaxA}{arg\,max}
\DeclareMathOperator*{\argminA}{arg\,min}
\DeclareMathOperator*{\maxA}{max}

\TheoremsNumberedThrough
\ECRepeatTheorems

\EquationsNumberedThrough

\MANUSCRIPTNO{OPRE-2022-08-445}

\begin{document}

\RUNAUTHOR{Bansak and Paulson}

\RUNTITLE{Dynamic Refugee Assignment}

\TITLE{Outcome-Driven Dynamic Refugee Assignment with Allocation Balancing}

\ARTICLEAUTHORS{
\AUTHOR{Kirk Bansak}
\AFF{Department of Political Science, University of California, Berkeley, \EMAIL{kbansak@berkeley.edu}} 
\AUTHOR{Elisabeth Paulson}
\AFF{Technology and Operations Management Unit, Harvard Business School,  \EMAIL{epaulson@hbs.edu}}
}

\ABSTRACT{
	This study proposes two new dynamic assignment algorithms to match refugees and asylum seekers to geographic localities within a host country. The first, currently implemented in a multi-year randomized control trial in Switzerland, seeks to maximize the average predicted employment level (or any measured outcome of interest) of refugees through a minimum-discord online assignment algorithm. The performance of this algorithm is tested on real refugee resettlement data from both the US and Switzerland, where we find that it is able to achieve near-optimal expected employment compared to the hindsight-optimal solution, and is able to improve upon the status quo procedure by 40-50\%. However, pure outcome maximization can result in a periodically imbalanced allocation to the localities over time, leading to implementation difficulties and an undesirable workflow for resettlement resources and agents. To address these problems, the second algorithm balances the goal of improving refugee outcomes with the desire for an even allocation over time. We find that this algorithm can achieve near-perfect balance over time with only a small loss in expected employment compared to the  employment-maximizing algorithm. In addition, the allocation balancing algorithm offers a number of ancillary benefits compared to pure outcome maximization, including robustness to unknown arrival flows and greater exploration.
}


\maketitle

\section{Introduction}

Host countries have, in recent years, been faced with increasing flows of refugees and asylum seekers. Currently, the United Nations Refugee Agency estimates that there are over 35 million refugees worldwide \citep{UNHCRtrends}. 
In most countries that accept refugees and/or asylum seekers, refugees and asylum seekers are assigned and relocated across various localities by migration authorities.
The capacities or target distributions of refugees across the localities are determined by authorities on a yearly or other regular basis. 

The goal of host countries is to help these new arrivals achieve economic self-sufficiency and other positive integration outcomes. Accordingly, a number of countries have begun to explore and implement outcome-based geographic matching in their refugee resettlement and/or asylum programs. Therefore, recent research studies the problem of efficiently assigning refugees to localities in order to maximize outcomes such as employment \citep{bansak2018improving, ahani2021placement}. This research falls within a broader area of policy interest as national resettlement programs seek new approaches to help ever-increasing flows of refugees and asylum seekers better integrate (e.g. find employment) in their host countries \citep[e.g.][]{mousa2018boosting, RePEc:hhs:lunewp:2018_007, golz2019migration, olberg2019enabling, acharya2020matching, ahani2021dynamic}. 

Outcome-based matching was introduced in the context of refugee and asylum-seeker assignment by \cite{bansak2018improving}, with the goal of leveraging administrative data to improve key refugee outcomes (e.g. employment in the host country) by optimizing refugees' geographic assignment within a country. To do so, machine learning methods are used to predict refugees' expected outcomes in each possible landing location as a function of the refugees' personal characteristics. Those expected outcomes are then used as inputs into constrained matching procedures to determine a location recommendation for each refugee. 

A greedy approach to the refugee assignment problem---one that assigns each refugee to the location (among those that are available) with the highest predicted outcome---is suboptimal when the resettlement locations have capacity constraints. This is the case in practice, where each location only has a certain number of slots in a given time period. For the United States, the time period is one year, but this can vary across host countries. Therefore, \cite{bansak2018improving} and other previous studies on outcome-based refugee matching \citep[e.g.][]{ahani2021placement,golz2019migration} have proposed optimal matching approaches to the refugee assignment problem that takes into account these capacity constraints.

This paper, along with the concurrent work of \cite{ahani2021dynamic}, are the first two papers to consider the \emph{dynamic} aspect of the outcome maximization matching problem. In many countries---including the United States (US), Switzerland, Sweden, the Netherlands, and Norway---refugees and asylum seekers must be assigned to a locality virtually immediately upon being processed by resettlement authorities. As a result, each arriving refugee or asylum seeker case (an individual or family) is typically assigned in an online fashion, and these assignments cannot be reversed. The dynamic aspect of this problem introduces a key trade-off between immediate and future rewards: assigning a current case to a location results in an immediate reward (namely, the employment score of the current case at that location), but also uses up a slot at that location for future arrivals.

This paper introduces two new dynamic matching algorithms. The first is a ``minimum-discord'' algorithm that seeks to maximize expected employment (or any alternative outcome of interest), and is currently employed in a pilot implementation in Switzerland, undertaken by the Swiss State Secretariat of Migration in collaboration with academic researchers. Details on the implementation in Switzerland are provided in Section \ref{sec:implementation}. The minimum-discord algorithm achieves near-optimal employment compared to the hindsight-optimal solution on real-world US and Swiss data. However, it can result in an imbalanced allocation to the localities over time which leads to implementation difficulties and an imbalanced workload for the resettlement offices.  

The second algorithm proposed in this paper is an extension that integrates principles of load balancing into the objective. Because each locality has a given amount of resources (e.g. resettlement officers and service providers) that cannot be transferred across localities, maintaining a steady workload is a first-order concern of resettlement agencies. Hence, building on the minimum-discord outcome maximization algorithm and borrowing ideas from queuing theory, the second algorithm incorporates wait time minimization into the assignment process. This allows refugees to be dynamically assigned to localities in a way that improves their expected employment scores while also maintaining a balanced allocation across the localities over time. Furthermore, the allocation balancing algorithm also offers ancillary benefits. In particular, it naturally handles the real-world scenario in which the total number of arrivals in a given period is not known in advance, and helps to improve the resilience of the underlying learning system through greater exploration.

This paper uses data from both the US and Swiss contexts to demonstrate the expected performance of the proposed approaches.  Furthermore, we discuss real-world constraints, phenomena, and difficulties that arose during Swiss implementation and our proposed solutions. 

\subsection{Contributions}

\begin{enumerate}
	\item \textbf{Minimum-discord outcome-maximizing dynamic assignment algorithm}. We propose a ``minimum-discord'' online algorithm that assigns arriving refugees to locations within a host country. The goal of the algorithm is to maximize the sum of individual outcomes along a horizon, while obeying the capacity constraints of each location. This is accomplished through a Monte-Carlo-sampling-based method that seeks to minimize the probability of choosing the ``wrong'' assignment in each time period compared to an offline benchmark. The proposed algorithm is a special case of the Bayes Selector algorithm \citep{vera2020bayesian}.
	
	\item \textbf{Allocation balancing dynamic assignment algorithm}. We demonstrate that an outcome-maximizing assignment (even a hypothetical implementation of the hindsight-optimal solution) can result in severe imbalance across the localities over time due to clustered arrivals of refugees with similar characteristics. Thus, we develop a second online algorithm that explicitly balances the trade-off between outcomes and having a balanced allocation to the localities over time using a single parameter, $\gamma$, that controls the weight placed on allocation balancing versus outcome maximization.
	 
	\item \textbf{Results on real refugee resettlement and asylum seeker data}. The results of the proposed methods are tested on real asylum seeker data from Switzerland and refugee resettlement data from one of the largest resettlement agencies in the US. In both cases, the proposed algorithms are able to improve upon the status quo assignment procedures by roughly 40-50\% and achieve 95-98\% of the hindsight-optimal solution.
	Using the allocation balancing algorithm, we demonstrate the trade-off between total employment and having a balanced allocation over time as $\gamma$ varies. In both contexts, we find that near perfect balance over time can be achieved with little loss in employment. 
	
	\item\textbf{Implementation details}. We describe practical constraints and learnings that arose during implementation in Switzerland. For example, we discuess how capacity updating throughout the year (resulting from uncertainty about the total number of individuals that will arrive each year) and a requirement to balance the geographic distribution of certain nationalities are treated in practice.
\end{enumerate}

\subsection{Related Literature} \label{sec:lit}

This paper is related to the existing literature on refugee assignment, online stochastic bipartite matching, and matching with queues. In what follows, we provide an overview of the most relevant literature from each stream.

\subsubsection{Geographic Assignment of Refugees and Asylum Seekers}

Prior research has proposed different schemes for refugee matching both across and within countries based on refugee and/or host location preferences \citep{fernandez2015tradable, moraga2014tradable,andersson2016assigning,delacretaz2016refugee,nguyen2021stability}. However, the lack of systematic data on preferences has thus far been a barrier to implementing these preference-based schemes.

In contrast, outcome-based matching was introduced in the context of refugee and asylum-seeker assignment by \cite{bansak2018improving}, with the goal of leveraging already existing data to improve key refugee outcomes (e.g. employment in the host country). However, the dynamic aspect of the problem is not considered by \cite{bansak2018improving}, nor by most previous studies on
outcome-based refugee matching \citep{ahani2021placement, golz2019migration, acharya2020matching}.
While \cite{RePEc:hhs:lunewp:2018_007} consider dynamically matching asylum seekers to localities, they focus on the goals of Pareto efficiency and envy-freeness across localities as opposed to outcome maximization.

\cite{ahani2021dynamic} is the closest to this paper. Like this paper, \cite{ahani2021dynamic} propose a dynamic matching algorithm to assign arriving refugees to locations within host countries with the goal of outcome maximization. The \emph{potentials method} proposed in \cite{ahani2021dynamic} is currently implemented by a resettlement agency in the US. For each newly arriving household, both the algorithm proposed in this paper and that of \cite{ahani2021dynamic} use a sampling procedure to solve many instances of the offline matching problem for the remaining horizon. \cite{ahani2021dynamic} then propose using dual variables from the offline problems to inform the assignment of the current arrival---a method referred to as the \emph{potentials} method. The algorithm proposed in this paper, on the other hand, assigns the current arrival to the location that minimizes the probability of a disagreement between the online algorithm and an offline benchmark. Both methods perform similarly on the data used in this paper. Our ``minimum-discord'' method, however, is both easily explainable and extends naturally to include allocation balancing, which is the focus of this paper.

Recent work also considers the relationship between the prediction and matching stages of dynamic refugee assignment \citep{kasy2023matching,bansak2023random} and group-fairness concerns \citep{freund2023group}.

\subsubsection{Stochastic Online Bipartite Matching}

Refugee matching is a special case of stochastic online bipartite matching, which has been a focus of operations and computer science researchers since the seminal work of \cite{karp1990optimal}. 

Two key features differentiate the refugee matching setting from the classic online matching problem. First, it is a weighted matching problem. Second, there is effectively an infinite number of arrival ``types,'' due to the large number of underlying covariates used to predict the outcome weights. While weighted online matching problems are well-studied, most existing methods rely on an assumption of finite types \citep{jaillet2012near,bumpensanti2020re,vee2010optimal,devanur2009adwords}. 
Although, in theory, the covariate domain could be discretized and adapted to a finite-type setting, this is undesirable. While there is prior research on distribution-free resource allocation problems, the performance guarantees of these algorithms nonetheless rely on a stationarity assumption \citep{devanur2019near}, which would not hold in practice in our setting. Rather, we seek to develop explainable methods that perform well, and do not focus on theoretical performance guarantees. 
The proposed method bears similarities to recent work by \cite{vera2020bayesian}.

\cite{vera2020bayesian} introduce a new framework for designing online policies given access to an offline benchmark. This framework is used to develop a meta-algorithm (``Bayes Selector'') for implementing low-regret online decisions across a broad class of allocation problems, including the assignment problem. In each state, the Bayes Selector chooses an action at each time interval that minimizes the likelihood of disagreement with an offline benchmark. 

When the number of arrival types is finite, \cite{vera2020bayesian} show that the Bayes Selector algorithm achieves constant regret for many special cases of the online assignment problem. This result is also proven in \cite{arlotto2019uniformly} for the multisecretary problem.
In this paper, we propose an outcome maximization algorithm that can be thought of as a special case of a Bayes Selector with infinite arrival types. When arrival types are drawn from a continuous distribution, \cite{bray2019does} shows that the multisecretary problem---which is a special case of the refugee matching problem with only two locations---no longer has bounded regret. Additionally, \cite{freund2019good} extend the methods introduced in \cite{vera2020bayesian} to more general decision-making problems, in particular showing that the uniform regret bound does not hold in settings with large uncertainty about the time horizon, which is likely to be the case in the refugee matching context.

\subsubsection{Allocation Balancing}

This paper develops an online matching algorithm that not only improves outcomes for refugees, but also balances the allocation to receiving locations (or, more generally, assignment options) over time. This aspect of the paper is related to one-sided matching with queues. In our setting, each location can be thought of as having a dedicated queue, since location assignments are made immediately and cannot be changed. 

A subset of online bipartite matching literature considers queuing systems. The topology of the queuing system is critical to the analysis method, and  most research in this area either focuses on optimally designing the underlying topology, or has topology that is substantially different from the refugee matching context (e.g., \cite{afeche2021optimal, leshno2019dynamic, vera2020dynamic}). 
 
\cite{balseiro2021regularized} propose an algorithm for online resource allocation that combines a welfare-maximizing objective with an arbitrary regularizer on the total consumption of each resource. This regularizer term can model what they call ``load balancing''---ensuring that the total level of consumption of each resource is balanced at the end of the horizon. While this has a similar flavor to our problem, we are interested in maintaining evenness in the allocation \emph{throughout} the horizon.

The kidney exchange literature also considers queueing models. For example, \cite{unver2010dynamic} develops an online mechanism for allocating kidneys with the goal of reducing wait time. \cite{bertsimas2013fairness} develop online kidney allocation policies that balance efficiency, fairness, and wait times. Recent work by \cite{ding2018fluid} also considers trade-offs between efficiency and fairness. However, unlike in our setting, the kidney exchange problem has a single queue. 

Because of the structure of the refugee matching problem (namely, the fact that each location has its own queue and decisions are irrevocable), the allocation balancing problem bears similarity to load balancing in computer science \citep{azar1998line}. 
However, the utility of load balancing algorithms is limited in our setting because of our additional goal of outcome maximization. Thus, in this paper, we develop a new approach that combines the objective of maximizing employment outcomes with achieving a balanced allocation over time.

The remainder of the paper is organized as follows. Section \ref{sec:setting_data} provides background on the refugee resettlement processes in the US and Switzerland and more details on the datasets used in this study. Section \ref{sec:problem} defines notation and describes the assumptions of the model and dynamics. Section \ref{sec:outcome_max} formulates the offline outcome maximization assignment problem, proposes an algorithm for the online setting, and demonstrates the performance of the method using the US and Swiss data. Section \ref{sec:implementation} provides further details of the implementation in Switzerland, including practical constraints and challenges. Section \ref{sec:balancing} introduces the allocation balancing component of the problem and proposes a new heuristic that balances employment outcomes and wait time. Section \ref{sec:conclusion} concludes.

\section{Settings and Data} \label{sec:setting_data}

This section provides more detail on the two specific contexts from which data are used in this paper: the refugee resettlement process in the US and the asylum procedure in Switzerland. The proposed methods are also applicable to many other countries where refugees and asylum seekers must be dynamically assigned to localities, including Sweden, the Netherlands, and Norway.

\subsection{Settings and Dynamics}

In the US context, we focus on the resettlement of UNHCR refugees, who are granted refugee status in the US prior to their arrival. In the US, the target number of refugees that will be resettled each year is determined by an annual cap set in advance of the start of the year. Refugees who are accepted into the US are then distributed across ten non-governmental resettlement agencies. Finally, each of those agencies maintains its own network of localities to which they assign newly arrived refugees, with capacities for each locality also determined in advance. 

In the Swiss context, we focus on asylum seekers, who request admission and asylum at a port of entry after entering a host country. In Switzerland, asylum seekers whose claims are not rejected are assigned on a case-by-case basis by the Swiss State Secretariat for Migration (SEM) to one of the 26 Swiss cantons. The assignment of asylum seekers across the cantons must follow an annually mandated proportionality key, which dictates the cantons' relative capacities to receive asylum seekers as a function of their population sizes. 

In both the US and Swiss contexts, the geographic placement for some refugees/asylum seekers is predetermined for reasons of family reunification, medical needs, or other special circumstances. For refugees and asylum seekers whose placement is not predetermined, decisions in both the US and Switzerland are driven primarily by capacity constraints at the locations, without a systematic attempt to optimize with respect to refugee/asylum-seeker outcomes. 
In Switzerland, the assignment to cantons is explicitly done on a quasi-random basis subject to the proportionality key.

Finally, the assignment batch size also varies by country. In Switzerland, the assignment is done on a one-by-one basis for each family after their post-arrival processing, and a number of other countries (e.g., the Netherlands) follow a similar procedure. In the US, assignment decisions are made on a weekly basis. Although the paper focuses on one-by-one assignment, Appendix \ref{sec:batching} discusses how the proposed methods can be readily extended to settings with batching. 

\subsection{Data and Scope}

For the US context, we use (de-identified) data on refugees of working age (ages 18 to 64) who were resettled in 2015-2016 into the US  by one of the largest US refugee resettlement agencies. For the Swiss setting, we use (de-identified) data on adult asylum seekers geographically assigned in 2015-2016 who eventually received full protection status specified under the Geneva Convention as well as those whose claim for Geneva protection status was rejected but were awarded subsidiary protection.

In both contexts, placement officers centrally assigned each case (individual or family) in the dataset to one of the possible locations---the 26 cantons in Switzerland, and about 30 resettlement locations in the US agency's network. Both datasets contain details on the refugees'/asylum seekers' characteristics (such as age, gender, origin, etc.), their assigned locations, and their employment outcome to be used for optimization. In the US context, the outcome is whether each refugee was employed 90 days after arrival at their assigned location. Refugees' employment status 90 days after their arrival is the key (and only) outcome metric that the resettlement agencies are required to report and that is tracked by the US government. In Europe, labor market integration is typically more challenging and takes longer for asylum seekers, and hence we use a longer-term employment outcome in the Swiss context. Specifically, we focus on whether or not each asylum seeker attained any employment within their first three years after assignment.

In both contexts, only ``free cases'' (those without prior family ties in their host country) are included in this study.  
This allows us to present a model and algorithm that aligns with the ongoing Swiss implementation which is scoped to include only free cases, as will be described in Section \ref{sec:implementation}. However, Appendix \ref{app:complexities} also shows how the proposed approaches can be extended to include cases with family ties, as may be the case in future implementations.

For each case, a vector of employment scores is constructed, where each element corresponds to the average probability for individuals within that case of finding employment (within 90 days for the US and within 3 years for Switzerland) if assigned to the particular location. To generate each case's outcome score vector, the same methodology is employed as in \cite{bansak2018improving}. Specifically, we use the data to generate models that predict the expected employment success of an individual at any of the locations, as a function of their background characteristics. These models were then applied to the cases who were assigned in 2015-2016 ($N=1,919$ for the US and $N=4,523$ for Switzerland) to generate their expected employment success at each location. This paper assumes that the employment scores are given for each case, and evaluates the proposed assignment algorithms relative to these predicted values.

The free cases that were assigned in 2016 ($N=1,175$ for the US and $N=1,502$ for Switzerland) are treated as the test cohorts in this paper. That is, the proposed algorithms are applied to these particular cohorts, in the specific order in which the families are logged as having actually arrived. The 2015 arrivals are utilized as historical data in the proposed algorithms. To further mimic the real-world process by which these cases would be assigned dynamically to locations, real-world capacity constraints are also employed such that each location can only receive the same number of cases that it actually received. 

\section{Notation and Model}\label{sec:problem}

Throughout, $[K]$ denotes the set of integers $\{1,...,K\}$, and $\mathbbm{1}\{\cdot\}$ denotes the indicator function. Additionally, $e_j$ denotes a vector with a value of one in the $j$-th component and zeros elsewhere. For a matrix $\mathbf{W}\in\mathbb{R}^{N_1\times N_2}$ and vector $\mathbf{w}\in\mathbb{R}^{N_2}$, $[\mathbf{w};\mathbf{W}]\in\mathbb{R}^{(N_1+1)\times N_2}$ denotes a new matrix whose first row is $\mathbf{w}$. 

We will assume throughout most of the paper that the total number of arrivals in a given year is known in advance. This assumption is generally not true, and is discussed further in Section \ref{sec:challenges}. In reality, the projected arrival numbers determined by resettlement authorities (for instance, the numbers projected by each of the ten US resettlement agencies in consultation with the US State Department) are revised throughout the year. Under this assumption, without loss of generality, we will assume that one case arrives each time period, and thus let $T$ denote both the number of arrivals and the time horizon. 
Let $M$ be the number of localities, indexed by $j$, with capacities/slots $s_j$, and $\sum_{j=1}^M s_j=T$. The capacities represent the number of individuals that each location can accommodate. 

The arriving cases are indexed by $t$. For simplicity of exposition, it will be assumed that each case is comprised of exactly one individual or, equivalently, that the capacities are set at the case-level (instead of the individual-level), which aligns with the ongoing implementation in Switzerland. This is further dicussed in Appendix \ref{sec:non-unit-case}, which also shows how the proposed methods can be extended to account for varying case sizes along with invidual-level capacities, as may be the situation in future implementations. We will let $a_{j}(t)$ be the number of cases allocated to location $j$ \emph{after} the allocation at time $t$, and define $\tilde{s}_j(t):=s_j-a_j(t)$ as the remaining slots at location $j$ after time $t$ (i.e., at the start of time $t+1$). 

The assignment of case $t$ to location $j$ results in a scalar outcome, $w_{tj}$. In the US context, the value of $w_{tj}$ represents the probability that case $t$ will find employment within 90 days if assigned to location $j$, and in the Swiss context it is the likelihood of finding employment within the first three years. In this paper, the outcome scores $w_{tj}$ are assumed to be known. In practice, they are estimated using a machine learning model that takes a large number of covariates as input \cite[see][]{bansak2018improving}.
In the online assignment problem, an arriving case is completely defined by its employment score vector, $\mathbf{w}_t$ (which is a function of the case's underlying covariates). Thus, we will use the matrix $\mathbf{W}$ with elements $w_{tj}$ to denote an arbitrary population of $T$ cases. Additionally, let $\mathbf{W}_t$ be shorthand for a population of  arrivals from time $t$ through $T$. We will assume that every free case can be assigned to any location with remaining capacity. In reality, even free cases may have idiosyncratic restrictions on which locations they can be assigned to (e.g., for medical reasons). This is further discussed in Appendix \ref{app:free-case-constraints}.

We will work in the underlying probability space $(\Omega, \mathcal{F}, \mathbb{P})$, where $\omega\in\Omega$ denotes a sample path of arrivals. Thus, there is a one-to-one correspondence between $\Omega$ and the set of all matrices $\mathbf{W}$, and fixing $\omega$ also fixes $\mathbf{w}_t$ for all $t\in[T]$. 
The vectors $\mathbf{\tilde{s}}(t-1)$ and set $\{\mathbf{w}_l\}_{l\in[t]}$ fully describe the state of the online assignment problem at time $t$. Therefore, let $S_t := (\mathbf{\tilde{s}}(t-1),\{\mathbf{w}_l\}_{l\in[t]})$ denote the state at time $t$. Note that if the arrivals in each time period are assumed to be independent, then the state could be described simply by $\mathbf{\tilde{s}}(t-1)$ and $\mathbf{w}_t$.  
To formalize the dynamics of the problem, the following features are assumed:\\
1. \textbf{Blind Sequentiality}: \emph{The cases are assigned in an order that is exogenously determined and unknown in advance, and each case $t$ must be assigned before case $t + 1$ is assigned.}\\
2. \textbf{Non-anticipativity}: \emph{Each case $t$ is assigned without knowledge of the outcome scores of the future arrivals.}\\
3. \textbf{Permanence}: \emph{Assignments cannot be changed once they are made.}

These features are representative of the real-world dynamics in many countries. Appendix \ref{sec:batching} demonstrates how batching, which violates the non-anticipativity assumption, can be incorporated into the proposed algorithms, resulting in performance gains. 

The binary variables $z_{tj}$ are the key decision variables, with $z_{tj}=1$ if case $t$ is assigned to location $j$ and $z_{tj}=0$ otherwise. Let $\Phi$ denote a full assignment of cases to locations such that the capacity constraints are satisfied, and let $\phi(t)$ denote the assignment for case $t$ (and thus $z_{t\phi(t)}=1$). Therefore, $w_{t\phi(t)}$ is the outcome of case $t$ under assignment $\Phi$, which could also be written as $\sum_{j\in[M]} z_{tj}w_{tj}$. The total employment score of matching $\Phi$ is given by
\begin{equation} \label{eq:totalscore}
w(\Phi) := \sum_{t=1}^T w_{t\phi(t)}=\sum_{t=1}^T \sum_{j=1}^M w_{tj} z_{tj}.
\end{equation}

\subsection{Queueing model} \label{sec:queuing}

To capture the allocation balancing problem, each location will be treated as a server with a dedicated queue. Although there are no physical queues, this modeling framework captures the relevant trade-offs. 
To that end, it is assumed that each location has a \emph{processing rate}, $\rho_j$, based on the resources (i.e. resettlement officers, service providers, and other related resources) at that location. This is the rate at which location $j$ can handle incoming cases. For example, if $\rho_j=1/2$, then location $j$ is able to handle one case every two periods on average.
Resettlement officers, service providers, and local community resources cannot be moved across locations. Therefore, for simplicity of expoisiton we assume that $\rho_j$ is stationary. However, it is straight-forward to adapt the analysis and proposed techniques to settings where $\rho_j$ varies over time or by features of the cases. We will assume throughout that capacities are set to be commensurate with processing rates, so that $\rho_j T = s_j$. Note that this assumption is essentially met by design in the resettlement program, as capacities for each location are programmatically decided on the basis of the resources at each location. However, in practice the value of $\rho_j$ could also be determined through interviews with case officers, particularly to understand case-level heterogeneities in processing rates.

The \emph{build-up} of location $j$ at time $t$, for $t> 2$, is given by 
\begin{equation} \label{eq:buildup}
b_{j}(t) = \max\{0, b_{j}(t-1)-\rho_j\} + z_{tj}
\end{equation}
with $b_j(1)=z_{1j}$ for all $j\in\{1,...,M\}$.
This is the build-up \emph{up to and including} the assignment at time $t$ but \emph{before} the processing at time $t$. This represents the number of cases either waiting or in process at time $t$. For each location, the ideal build-up level is in the interval $(0,1]$, indicating that the location is actively settling a case and no cases are waiting. When $b_j(t)>1$, cases are waiting to be processed at location $j$, and when $b_{j}(t)=0$ location $j$ is idle. 

\section{Outcome Maximization}\label{sec:outcome_max}

This section proposes a minimum-discord online assignment algorithm that seeks to maximize the sum of outcome scores across the horizon. In this section, the build-up at each location is not considered. Section \ref{sec:balancing} will extend this algorithm by proposing a modified version that additionally seeks to minimize build-up.

First, we introduce the offline version of the outcome maximization problem. For a given set of arrivals $\mathbf{W}$, the offline optimization problem is:
\OneAndAHalfSpacedXI
\begin{equation} \label{prob:offline} \tag{\textsc{OutcomeMax}}
	\begin{aligned}
		\maxA_{\mathbf{Z}} &\sum_{t=1}^{T}\sum_{j=1}^M w_{tj} z_{tj}\\
		s.t. 	& \sum_{j=1}^M z_{tj} = 1 \;\;\; \forall \; t \in [T]\\
		& \sum_{t=1}^T z_{tj} = s_j \;\;\; \forall \; j \in [M]\\
		& \mathbf{Z}\in\{0,1\}^{T\times M}
	\end{aligned}
\end{equation}
\DoubleSpacedXI
where $\mathbf{Z}$ is the assignment matrix with elements $z_{tj}$.
The solution to \ref{prob:offline} is the outcome maximizing assignment for a population $\mathbf{W}$. When a particular population or sample path is specified, we may write this problem as \ref{prob:offline}$(\mathbf{W})$ or \ref{prob:offline}$(\omega)$, and its optimal objective value represents an upper-bound for any assignment of that population or sample path. It is well known that an optimal solution to \ref{prob:offline} can be found by solving the linear programming (LP) relaxation of \ref{prob:offline} \citep{bertsimas1997introduction}. Thus, solving \ref{prob:offline} is generally fast (e.g., for $T\leq 3,000$ \ref{prob:offline} can be solved in less than one second). See \ref{app:scale} for detailed run-time metrics.

The true online assignment problem is a dynamic program. In other words, the algorithm must make an assignment, given the current state, without knowledge of the outcome score vectors of future arrivals.
Because of the online nature of the problem, it is helpful to let the notation $\ref{prob:offline}(\mathbf{W}_t, \tilde{\mathbf{s}}(t-1))$ describe solving \ref{prob:offline} for time steps $t$ onward for population $\mathbf{W}_t$, starting with capacities $\tilde{\mathbf{s}}(t-1)$.

In theory, the optimal solution to the dynamic problem could be found by solving Bellman's equation, given by

\begin{equation} \label{prob:bellman}
	\begin{aligned}
	V_t(S_t)=\max_{\phi(t)\in[M]} &\left(w_{t\phi(t)} + \int_{\omega \in \Omega} \mathbb{P}(\omega|S_t)V_{t+1}(\mathbf{\tilde{s}}(t-1)-e_{\phi(t)}, \{\mathbf{w}_l\}_{l\in[t]}\cup \mathbf{w}(\omega)_{t+1}) \right)\\
	s.t.& \; e_{\phi(t)}\leq \tilde{s}_j(t-1) \;\;\forall \; j\in[M]
	\end{aligned}
\end{equation}
The optimal policy is the maximizer of the right-hand side of the equation above. Due to the so-called ``curse of dimensionality'' \citep{bellman1966dynamic} arising from the large number of locations and continuous outcome scores, Problem \ref{prob:bellman} cannot be solved directly, even if the probabilities $\mathbb{P}(\omega|S_t)$ were known.
Many heuristics and approximation methods have been proposed to solve Problem \ref{prob:bellman}. Our chosen solution method, a special case of the Bayes Selector method introduced in \cite{vera2020bayesian}, is described in the following section.

\subsection{Minimum-Discord Online Algorithm}

Let
\begin{equation}
Q(\phi(t),S_t):=\{\omega\in\Omega: \phi(t)\notin  \argmaxA_{j}  \left(w_{tj} + V_{t+1}(\mathbf{\tilde{s}}(t-1)-e_{j}, \{\mathbf{w}_l\}_{l\in[t]}\cup \mathbf{w}(\omega)_{t+1}) \right)\}
\end{equation}  
be the event that assigning case $t$ to location $\phi(t)$ is not optimal according to \ref{prob:offline}$(\omega)$. This definition allows for the possibility that there are multiple optimal decisions according to the offline benchmark.
Furthermore, let
$q(\phi(t),S_t):=\mathbb{P}[Q(\phi(t),S_t)|S_t]$
be the \emph{disagreement probability}. The most general version of the Bayes Selector algorithm proposed by \cite{vera2020bayesian} chooses the location at time $t$ that minimizes $q(\phi(t),S_t)$.
The algorithm proposed in this paper chooses the location which minimizes an approximation of these disagreement probabilities in each time period. This approach is referred to as \emph{minimum-discord}, since the goal is to minimize the likelihood of disagreement with the offline optimal solution at time $t$. We note that this method does not take into account the degree of disagreement. An alternative approach could select the location that minimizes the expected optimality gap, as opposed to minimizing the likelihood of making a suboptimal decision. This is elaborated on in Appendix \ref{sec:min-risk}. 

\cite{vera2020bayesian} establish performance guarantees for the Bayes Selector algorithm in many settings; however, the assumptions that underlie these guarantees do not hold in our setting which places no assumptions on the underlying arrival distribution. The focus of this paper is on proposing explainable algorithms with strong empirical performance on the real-world setting. Nonetheless, in \ref{app:regret} we provide a characterization of the expected regret of any online algorithm in terms of the disagreement probabilities, following Lemma 1 of \cite{vera2020bayesian}.

Because there are no assumptions placed on the arrival process, we use a Monte Carlo sampling procedure to estimate $q(\phi(t),S_t)$. The intuition is as follows.  When case $t$ arrives, we generate $K$ random trajectories of future arrivals $t+1$ through $T$, denoted by $\{\mathbf{W}_{t+1}^k\}_{k=1}^K$. For each random trajectory $k \in[K]$, the offline problem \mbox{\ref{prob:offline}$_t([\mathbf{w}_t;\mathbf{W}_{t+1}^k],\tilde{\mathbf{s}}(t-1))$} is solved. 

Let $n_j(t)$ be the number of times that case $t$ is assigned to location $j$ across the $K$ trajectories. The quantity $1-q(j,S_t)$---namely, the probability that location $j$ is an optimal action---is approximated by $n_j(t)/K$. Therefore, minimizing our approximation of $q(j,S_t)$ is equivalent to assigning case $t$ to location $\argmaxA_j n_j(t)$, that is, the location that they were assigned to most often in the random instances.
The proposed method is formally defined below.
\begin{heuristic}[\textsc{MinDiscord}] \label{heur:modal}
	Case $t$ is assigned to location 
	$$\phi(t) := \argmaxA_{j\in[M]} \sum_{k=1}^Kz^{k}_{tj},$$ with ties broken randomly, where  
	\begin{equation*} 
			\mathbf{Z}^{k} = \argmaxA_{\mathbf{Z}} \text{\ref{prob:offline}}([\mathbf{w}_t;W^k_{t+1}],\tilde{\mathbf{s}}(t-1)).
	\end{equation*}
\end{heuristic}

The Monte Carlo sampling approach requires a ``sampling population'' to draw sample from, which we denote by $\mathcal{A}$. In this paper, $\mathcal{A}$ is comprised of the 2015 arrivals. Algorithm \ref{alg:modal}, defined below, is the online assignment algorithm that employs Method \ref{heur:modal} in each time period. We note that the choice of $\mathcal{A}$ should depend on the level of non-stationarity in the arrival process. If the data is highly non-stationary, $\mathcal{A}$ could be comprised of a shorter, more recent window of arrivals.

\begin{varalgorithm}{\textsc{OnlineMinDiscord}}
	\SingleSpacedXI
	\caption{(Min-Discord Online Assignment)}\label{alg:modal}
	\begin{algorithmic}[1]
		\State \textbf{initialize} $\tilde{s}_j(0)\leftarrow s_j$ for all $j\in\{1,...,M\}$
		\For{$t$ in $1,...,T$}
		\For{$k$ in $1,...,K$}
		\State $\mathbf{W}_{t+1}^k \leftarrow$ $T - t$ randomly drawn cases from set $\mathcal{A}$ with replacement
		\State $\mathbf{Z}^k \leftarrow \argmaxA$  \ref{prob:offline}($[\mathbf{w}_t;\mathbf{W}_{t+1}^k],\tilde{\mathbf{s}}(t-1)$)
		\EndFor
		\State $\phi(t) \leftarrow\argmaxA_j \sum_{k}z^k_{tj}$ (with ties broken randomly)
		\State $\mathbf{\tilde{s}}(t) \leftarrow\mathbf{\tilde{s}}(t-1)-e_{\phi(t) }$
		\EndFor\\
		\Return $\Phi^{MD}=\{\phi(t)\}_{t=1...T}$
	\end{algorithmic}
\end{varalgorithm}

\subsection{Performance of \ref{alg:modal}}

\begin{figure}[ht!]
	\centering
	\caption{Results of online algorithms for US refugees (left) and Swiss asylum seekers (right) in 2016.} \label{fig:results}
	\includegraphics[width=.49\textwidth]{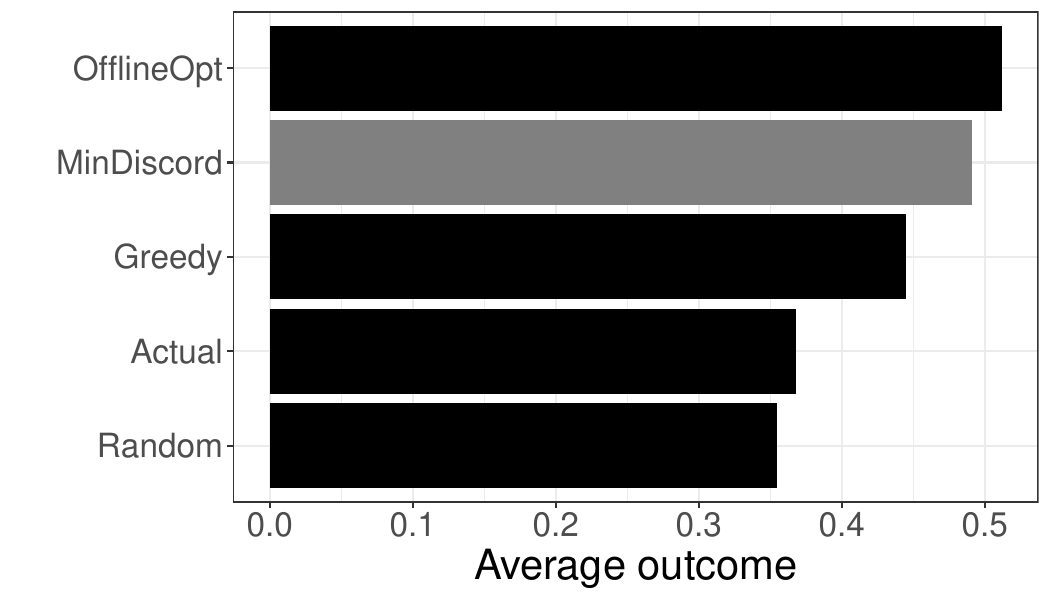}
	\includegraphics[width=.49\textwidth]{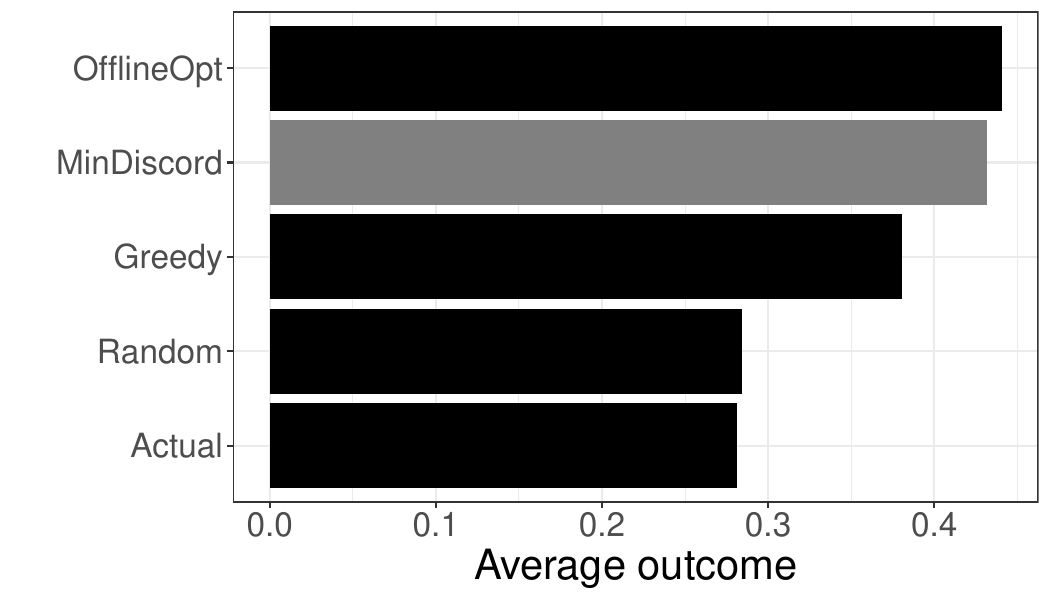}
\end{figure}

Figure \ref{fig:results} shows the results of applying \ref{alg:modal} to the 2016 arrivals (both US and Swiss). Throughout the paper, unless otherwise specified, we use $K=5$ for \ref{alg:modal}. We compare \ref{alg:modal} to four benchmarks: the actual historical assignment, the hindsight-optimal solution, greedy assignment, and random assignment. The first benchmark assigns each case to the location that they were assigned to in reality under the status quo procedures. Although, for this benchmark, we could measure employment according to whether or not the cases actually found employment in reality (since this is contained in the data), for all benchmarks we measure employment according to the predicted employment scores, $\mathbf{W}$, so that they are all evaluated with respect to the same metric. (We note, however, that using actual employment results in an almost identical total employment score for this benchmark.)

The hindsight-optimal solution, OfflineOpt, is included as a benchmark because, while it cannot be performed in a real-world dynamic context, it sets an upper bound of what is achievable by any algorithm. In the greedy algorithm, each case is assigned sequentially to the location with the highest expected employment score for that case, out of locations with remaining capacity. 
Finally, the employment score under random assignment for case $t$ is given by $\sum_{j\in[M]} w_{tj} \frac{s_j}{T}$, which we include as a simple reference point. A comparison of \ref{alg:modal} to the method proposed by \cite{ahani2021dynamic} is also included Appendix \ref{sec:additional_computations}, though we note that the methods perform quite similarly.

Figure \ref{fig:results} shows the results. On the US data, \ref{alg:modal} achieves 96\% of the employment score of the hindsight-optimal solution. This is compared to the greedy, random, and actual historical assignment benchmarks, which achieve 87\%, 69\%, and 72\% of the hindsight optimal employment levels, respectively. 
On the Swiss data, \ref{alg:modal} achieves 98\% of the employment score of the hindsight-optimal solution. In this case, the greedy, random, and actual historical assignments achieve 86\%, 64\%, and 64\% of the hindsight optimal solution, respectively. Outcomes by certain subgroups (e.g., nationality and sex) are shown in Figure \ref{fig:subgroups}.

The optimality gap of \ref{alg:modal} is primarily because of nonstationarity in the arrival processes. When the arrival dates of the cases are randomly perturbed and $\mathcal{A}$ is taken to be the 2016 arrivals---mimicking a stationary process---the optimality percentage of the proposed algorithm compared to the hindsight-optimal solution increases to about 99.5\% on both the US and Swiss data. This was calculated as the average optimality percentage across fifty random instances, where in each instance the arrival dates of the cases are randomly shuffled; in each of these instances, the optimality percentage was between 99.4\% and 99.7\%. However, the focus of this paper is on the performance of the proposed algorithms on the real, non-stationary, arrival data.

\section{Implementation Details} \label{sec:implementation} 

This section provides details on the current pilot implementation of \ref{alg:modal} in Switzerland. Additionally, we discuss implementation complexities that motivated the development of a second algorithm, described in Section \ref{sec:balancing}.

\subsection{Background}

In coordination with the SEM in Switzerland and a multi-university collaboration between researchers from ETH Zurich, Stanford, Dartmouth, Harvard, and the University of California, Berkeley, a multi-year pilot implementation of \ref{alg:modal} is ongoing in Switzerland. The pilot began in January 2020 and is projected to end in 2024. As described further below, the pilot includes a randomized control trial (RCT) and targets the optimization of three-year employment outcomes; for this reason, results are not yet available, and the final results will not be available until three years after the completion of the pilot. The objective of the pilot is to generate rigorous evidence of impact on asylum-seeker employment, based upon which a broader and more permanent implementation of these methods can then be considered by the SEM.

The pilot implementation applies to all adult asylum seekers (or families that include at least one adult) who (a) obtained subsidiary or Geneva convention protection status (and hence who are granted asylum and allowed to stay in Switzerland), (b) who are free to be assigned to any canton (i.e., do not have pre-existing family ties, medical constraints, or other special arrangements), and (c) are part of the ``accelerated procedure'' track in the Swiss asylum process. The accelerated procedure is used for relatively uncomplicated cases whose status---whether they will be granted asylum or will be removed from Switzerland---can be designated in a relatively prompt manner, with a target of less than 100 days. 

\subsection{Pilot Set-up}

As described earlier, placement officers in Switzerland are in charge of determining the cantonal assignment of asylum seekers. In our pilot implementation, the placement officers have been provided with specialized software that generates a recommended canton for each asylum seeker case (i.e. family or individual). The placement officers maintain the ability to override the recommendation if necessary, but they are encouraged to take the recommendation; as mentioned, the pilot scope includes only asylum seeker cases that can be assigned to any canton. The RCT design is simple: each asylum seeker case that will be assigned is first randomly allocated to either the control or treatment condition. In the control condition, the canton recommendation is generated randomly. In the treatment condition, the canton recommendation is generated via \ref{alg:modal}.

The distribution of asylum seekers in Switzerland follows a cantonal proportional distribution key. Accordingly, in our implementation, the assignment of asylum seeker cases is subject to canton capacity constraints that follow this proportional distribution key, which is enforced separately for the treatment and control cases. In other words, the treatment and control cases have fully independent capacity at each canton to limit interference in the RCT. The capacity for cases that are out of scope (e.g., cases with family ties) are also independent from the pilot. Furthermore, capacities are set for each of the treatment and control cohorts at the case-level in the pilot implementation. Hence, the implementation of \ref{alg:modal} is applied as described in the text with cases as the units of interest, though with one additional consideration: the proportional distribution constraints must be achieved independently for six different nationality groups in accordance with a Swiss legal requirement. The six groups are comprised of asylum seekers from: (1) Afghanistan, (2) Turkey, (3) Georgia, (4) the Maghreb countries, (5) a handful of specially identified countries (Albania, Benin, Burkina Faso, Bosnia and Herzegovina, Ghana, Guinea, Gambia, India, Moldova, the Republic of North Macedonia, Mongolia, Nigeria, Kosovo, Senegal, and Serbia), and (6) all other countries. The need to achieve proportional distribution independently for each of these six groups is tantamount to, and hence is achieved by, implementing \ref{alg:modal} separately and independently for each of these six groups.

\subsection{Challenges}\label{sec:challenges}

Because there is no advance processing prior to the arrival of asylum seekers, flows of asylum seekers can be somewhat unpredictable. Regional and global events, such as conflicts and wars, can lead to sudden changes in the types of asylum seekers who are arriving and their rate of arrival. Relative to the context of assigning UNHCR resettled refugees, this poses a more significant challenge for setting and controlling the capacity constraints, given that the number of cases that will need to be assigned by the end of the year (or within any period of time) is fundamentally uncertain. This uncertainty results in a violation of the modeling assumption that $T$---the total number of arrivals---is known in advance.

Nonetheless, ensuring that the distribution of the assignments across cantons meets the proportional allocation key by the end of each calendar year is of critical administrative importance, which requires that our capacity targets not exceed the actual number of annual arrivals without knowing what that number will be in advance. We employ intermittent updating of the capacity constraints to deal with this challenge in the pilot implementation. Because resources to process and receive new asylum seekers within each canton are limited and cannot freely move across cantons, it is also important that assignments to any given canton are not too concentrated within a period of time (e.g. if a canton's quota for the entire year were assigned to it in a single month).

To deal with both of these issues, we employ recent trends to project the number of arrivals in shorter intervals (e.g. 1-4 months) and intermittently add capacity according to the proportional allocation key over the course of the year. In doing so, we are able to avoid overloading any canton and protect against a divergence from the proportional allocation key. The cost, however, is inefficiency in two regards. First, the updating process itself entails analyses that cannot be easily automated, and hence requires additional human labor. Second, introducing smaller chunks of capacity intermittently over time can cut into the ability of the algorithm to maximize gains. 

These considerations and learnings from the pilot implementation have thus motivated our proposal for a second algorithm (\ref{alg:modal_balance}), presented in the following sections, that maintains a balanced geographic distribution over time. By incorporating this allocation balancing component, \ref{alg:modal_balance} not only ensures that all locations have a steady stream of arrivals throughout the year, thus ensuring that local resources in any location are not outstripped by a sudden imbalanced influx of arrivals at any given time, but also naturally ensures that the overall distribution of cases will meet the proportionality targets regardless of uncertainty in the arrival numbers.

\section{Allocation Balancing}\label{sec:balancing}

Motivated by learnings from the pilot implementation of \ref{alg:modal}, this section presents an extension that strives to maintain a balanced, proportional allocation over time to each locality by considering each locality to be a server with a dedicated queue (see Section \ref{sec:queuing} for the modeling details).

\subsection{Imbalance Under Outcome Maximization}\label{sec:resulting_imbalance}

\begin{figure}[h!]
	\begin{center}
		\caption{Allocation to nine largest locations over time for US data (left) and Swiss data (right) using \ref{alg:modal}.} \label{fig:usresults_alg_imbalance}
		\includegraphics[width=.49\textwidth]{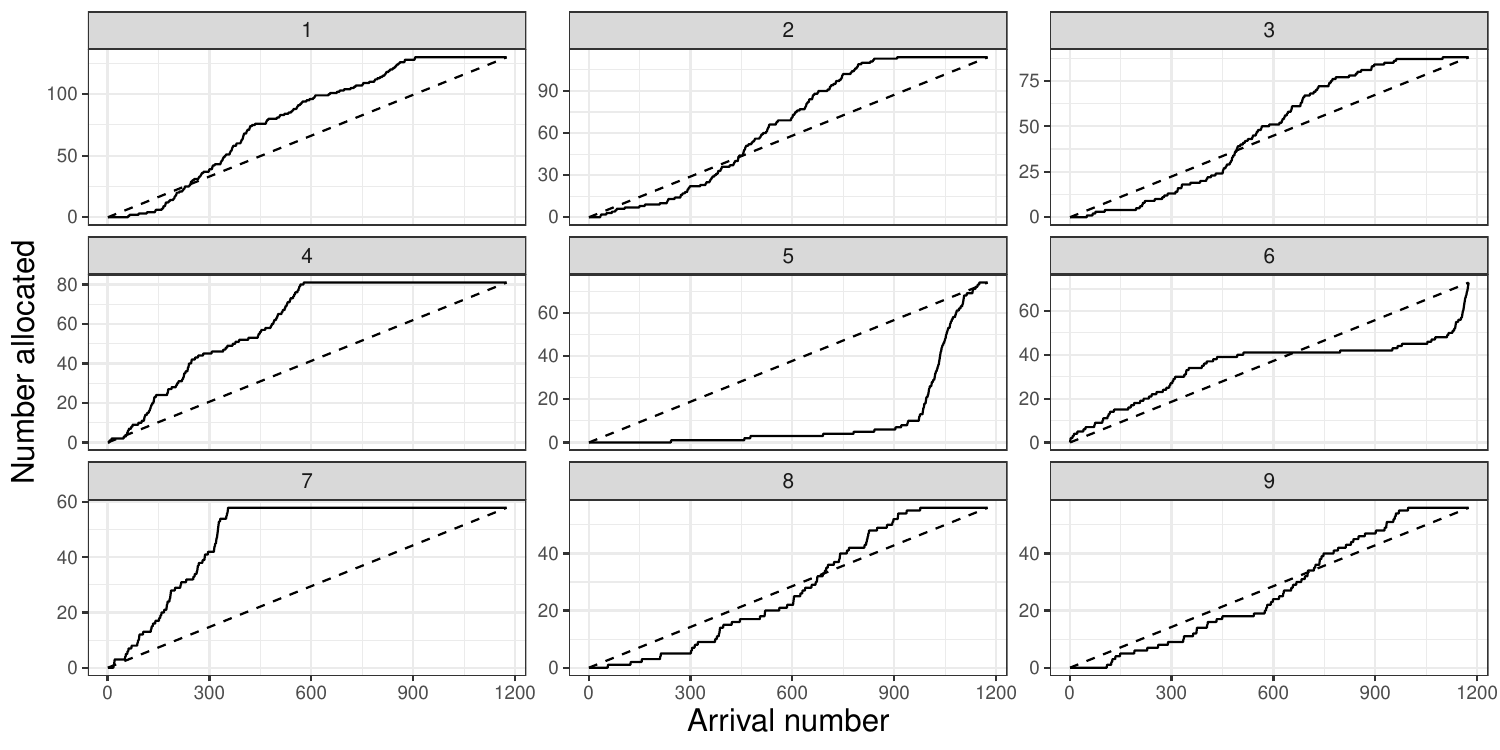}
		\includegraphics[width=.49\textwidth]{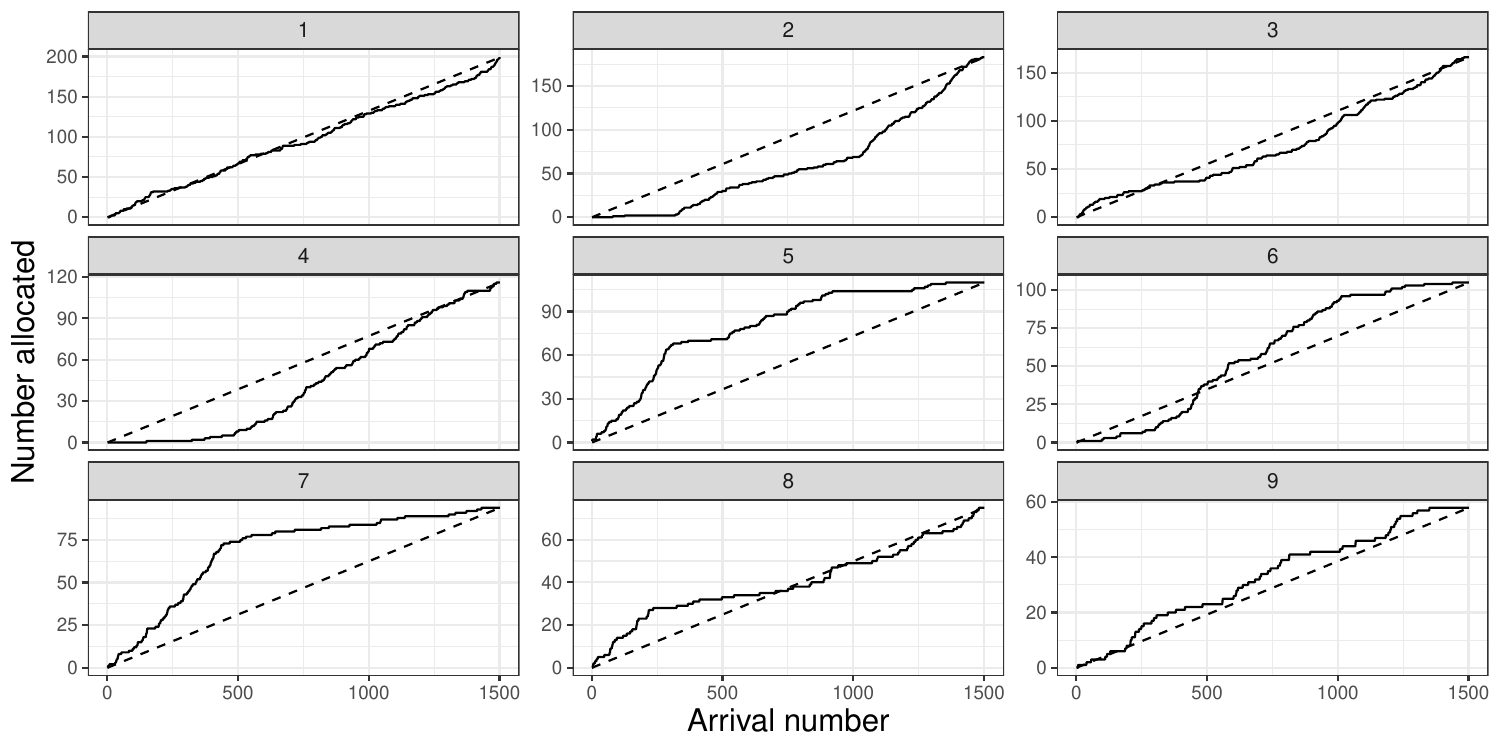}
	\end{center}
\end{figure}

Although \ref{alg:modal} performs well in terms of maximizing outcomes, it results in significant imbalance in the allocation to localities over time. Figure \ref{fig:usresults_alg_imbalance} shows the cumulative allocation to the largest nine locations over the horizon for the US data (left) and Swiss data (right) obtained by using \ref{alg:modal}. 

For the US data, the average queue length (defined as $\max(b_j(t)-1,0)$) across all locations resulting from \ref{alg:modal} is 7.5. For comparison, the average queue length under the actual historical assignment is 3.1. Thus, switching to an optimization approach does indeed lead to longer queues and wait times than under the status quo procedure.  
Similarly, for the Swiss data, the average queue length of \ref{alg:modal} is 6.3 compared to 3.3 under the actual historical assignments.
We note that this is not simply a consequence of the particular choice of online algorithm, nor entirely a consequence of the online nature of the problem: even the hindsight-optimal solution results in imbalance over time (see Figure \ref{fig:imbalance_opt}). 

This imbalance is primarily driven by non-stationarity in the arrival process. Indeed, when \ref{alg:modal} is applied to the same 2016 data but with a randomly perturbed arrival sequence (and with $\mathcal{A}$ set to be the 2016 arrival cohort), mimicking a stationary process, the average queue length across five random instances is 1.7 (see Section \ref{sec:additional_computations}). Because refugee inflows are, in part, due to international events, there can be clustering of arrivals with specific background characteristics---particularly with respect to country of origin, which is one of the predictors that underlies the employment scores. This can lead to clustering in the subsequent assignment, causing imbalance. As described in Section \ref{sec:implementation}, this phenomenon encourages a conservative capacity updating approach in order to ensure that no location exceeds their proportionality key at the end of the horizon. An imbalanced allocation is also highly undesirable for resettlement service providers who cannot move between locations. The allocation balancing method, described in the following sections, mitigates these issues and provides ancillary benefits.

\subsection{Offline Benchmark}

Because the number of slots at each location is fixed, minimizing queue length/wait time is effectively equivalent to minimizing wait time \emph{and} idle time. Therefore, we focus on minimizing wait time explicitly, while also noting the subsequent impact of the proposed methods on idle time. 
For simplicity of exposition, we will assume that the cost of wait time is identical across locations, although extending the algorithm to the non-identical case is straightforward. 

First, consider a new variant of the offline benchmark that penalizes wait time, given by:

\OneAndAHalfSpacedXI
\begin{equation} \label{prob:offline_balance} \tag{\textsc{Balance}}
	\begin{aligned}
	\maxA_{\mathbf{Z},\mathbf{b}} &\sum_{t=1}^{T}\sum_{j=1}^M w_{tj}z_{tj} - \gamma  \sum_{t=1}^{T}\sum_{j=1}^M \lceil b_{j}(t)-1\rceil \mathbbm{1}\{b_{j}(t)> 1\} \\
	\text{s.t.} 	& \sum_{j\in[M]} z_{tj} = 1 \;\;\; \forall \;\;\; t \in [T]\\
	& \sum_{t\in[T]} z_{tj} = s_j\;\;\; \forall \;\;\; j \in [M]\\
	& b_{j}(t) = \max\{0,b_{j}(t-1)-\rho_j\} + z_{tj} \;\;\; \forall \;\;\; t \in \{2,...,T\},\; j\in[M]\\
	& b_j(1)=z_{1j} \;\;\; \forall \;\;\; j \in [M]\\
	& \mathbf{Z}\in\{0,1\}^{T\times M}
	\end{aligned}
\end{equation}

\DoubleSpacedXI
Recall that $b_j(t)$ denotes the build-up at location $j$ at time $t$, and $\rho_j$ is again the processing rate of location $j$.
In the objective function of \ref{prob:offline_balance}, wait time cost is incurred when $b_{j}(t)> 1$, and $\lceil b_j(t)-1\rceil$ is the number of cases waiting at time $t$.
The parameter $\gamma$ (assumed to be non-negative) is a weight that balances the trade-off between outcomes and wait time, and can be thought of as the relative cost of wait time. In practice, this parameter could be set either according to a cost-benefit analysis such that the units of measure were commensurate with one another, or according to an empirically driven decision on a value that results in acceptable balance across locations over time. 

Let \ref{prob:offline_balance}$(\mathbf{W}_{t},\tilde{\mathbf{s}}(t-1), \mathbf{b}(t-1))$ denote solving \ref{prob:offline_balance} from time $t$ onward, for population $\mathbf{W}_t$ with capacities $\tilde{\mathbf{s}}(t-1)$ and initial buildup $\mathbf{b}(t-1)$.
Recall that in \ref{alg:modal}, \ref{prob:offline} is solved $K$ times for each new arrival, each time using a randomly generated sample of future arrivals. This same approach will be used to develop the new online allocation balancing assignment algorithm.

However, unlike \ref{prob:offline}, \ref{prob:offline_balance} cannot be solved to optimality as a linear program. The variables $b_j(t)$ are defined by non-linear expressions, the objective function of \ref{prob:offline_balance} is non-linear, and finally the assignment variables are binary. Due to advances in mixed-integer programming (MIP), \ref{prob:offline_balance} can still be solved using state-of-the-art MIP solvers, and one can obtain partial speed-ups by linearizing and relaxing parts of the problem. However, these approaches nonetheless result in substantially increased run-time compared to \ref{prob:offline} (see \ref{app:scale} for further discussion).
Thus, instead of using \ref{prob:offline_balance} as our offline problem, we propose an alternative method that uses a greedy version of \ref{prob:offline_balance}. We show that this approach results in strong empirical performance and argue why a greedy approach is reasonable for allocation balancing.

\subsection{Online Allocation Balancing Algorithm}

In this section we propose a greedy version of \ref{prob:offline_balance} to use as the offline problem in the online allocation balancing algorithm. In an online setting, the past assignments to each location are readily observable. Thus, at time $t$, the online algorithm has access to $b_j(t-1)$ for all locations $j$. Consider the following problem at time $t$:

\OneAndAHalfSpacedXI
\begin{equation} \label{prob:offline_balance_greedy} \tag{\textsc{GBalance}}
	\begin{aligned}
		\maxA_{\mathbf{Z}} &\sum_{l=t}^{T}\sum_{j=1}^M w_{lj}z_{lj} -  \gamma\sum_{j=1}^M z_{tj}  \left\lceil\frac{b_j(t-1)-\rho_j}{\rho_j}\right\rceil\mathbbm{1}\{b_j(t-1)>0\} \\
		\text{s.t.} 	& \sum_{j\in[M]} z_{lj} = 1 \;\;\; \;\;\;\forall \;\;\; l \in \{t,...,T\}\\
		& \sum_{l=t}^T z_{lj} = \tilde{s}_j(t) \;\;\; \forall \;\;\; j \in [M]\\
		& \mathbf{Z}\in\{0,1\}^{N\times M}
	\end{aligned}
\end{equation}
\DoubleSpacedXI
\noindent \ref{prob:offline_balance_greedy} takes $\mathbf{b}(t-1)$ as input, and weights the employment score of case $t$ by the wait time cost incurred by case $t$. The wait time that case $t$ experiences if assigned to location $j$ is the length of time until all earlier cases are done being processed, starting from time $t$, namely $\left\lceil\frac{b_j(t-1)-\rho_j}{\rho_j}\right\rceil$. 
Because $\mathbf{b}(t-1)$ is known prior to the $t$-th arrival, \ref{prob:offline_balance_greedy} has a linear objective function. 
Thus, as with \ref{prob:offline}, the optimal solution to \ref{prob:offline_balance_greedy} can be found by solving its LP relaxation. In fact, solving \ref{prob:offline_balance_greedy} is as fast as solving \ref{prob:offline}, making this problem appealing for use in an online setting. 

 To build intuition for  \ref{prob:offline_balance_greedy}, we present the following lemma, which bridges \ref{prob:offline_balance} and \ref{prob:offline_balance_greedy}. The proof of Lemma \ref{lemma:equivalence} can be found in \ref{app:lemma-proof}.

\begin{lemma}\label{lemma:equivalence}
	The objective function of \ref{prob:offline_balance} is equivalent to
	\begin{align}\label{exp:offline_balance_mod}
		\sum_{t=1}^T\sum_{j=1}^M w_{tj}z_{tj}-\gamma \sum_{t=1}^T \sum_{j=1}^M  z_{tj} \left\lceil\frac{b_j(t-1)-\rho_j}{\rho_j}\right\rceil\mathbbm{1}\{b_{j}(t-1)>0\}.
	\end{align}
\end{lemma}

Notice that the objective function of \ref{prob:offline_balance_greedy} is a \emph{greedy} version of Expression \ref{exp:offline_balance_mod}---namely, it does not calculate wait time for the entire horizon, but does so only for the current arrival (hence the name \textsc{\textbf{G}reedy} \textbf{\textsc{Balance}}). 

Although the greedy method does not work well when it comes to outcome maximization, it does work well for minimizing wait time.
In terms of wait time, taking a slot from location $j$ immediately increases the build-up at location $j$. This makes it less likely for arrivals in the near-future to be assigned to location $j$, which could be consequential especially if, for some of these arrivals, location $j$ is highly desirable. However, this effect is short-lived: it is only relevant if an arrival in the near future (i.e., before the current arrival can be fully processed) would also be assigned to location $j$. Thus, if there are many locations, or if location $j$ has few slots (implying that the probability of any given arrival being assigned to location $j$ is small), this effect is mitigated.

Accordingly, we propose a new assignment method. Method \ref{heur:modal-balancing} is similar to Method \ref{heur:modal}, but assigns the current arrival based on the solution to \ref{prob:offline_balance_greedy} instead of \ref{prob:offline} as in Method \ref{heur:modal}.

\begin{heuristic}[Allocation-Balancing MinDiscord] \label{heur:modal-balancing}
	Case $t$ is assigned to location 
	$$\phi(t) = \argmaxA_{j\in[M]} \sum_{k\in[K]}z^{k}_{tj},$$ with ties broken randomly, where 
	\begin{equation*}
		\begin{aligned} 
			\mathbf{Z}^{k} = \argmaxA_{\mathbf{Z}} \; &\text{\ref{prob:offline_balance_greedy}}\left([\mathbf{w}_t;\mathbf{W}^k_{t+1}],\tilde{\mathbf{s}}(t-1), \mathbf{b}(t-1)\right).
		\end{aligned}
	\end{equation*}
\end{heuristic}

The online algorithm based on Method \ref{heur:modal-balancing} is presented as Algorithm \ref{alg:modal_balance} below.

\begin{varalgorithm}{\textsc{OnlineBalance}}
	\SingleSpacedXI
	\caption{(Allocation-Balancing Online Assignment)}\label{alg:modal_balance}
	\begin{algorithmic}[1]
		\State \textbf{initialize} $\tilde{s}_j(0)\leftarrow s_j$ for all $j\in\{1,...,M\}$
		\For{$t$ in $1,...,T$}
		\For{$k$ in $1,...,K$}
		\State $\mathbf{W}_{t+1}^k \leftarrow$ $T - t$ randomly drawn cases from set $\mathcal{A}$
		\State $\mathbf{Z}^k \leftarrow \argmaxA$ \ref{prob:offline_balance_greedy}$([\mathbf{w}_t;\mathbf{W}^k_{t+1}],\tilde{\mathbf{s}}(t-1), \mathbf{b}(t-1))$
		\EndFor
		\State $\phi(t)\leftarrow\argmaxA_j \sum_{k\in[K]}z^k_{tj}$ (ties broken randomly)
		\State $\tilde{\mathbf{s}}(t) \leftarrow \tilde{\mathbf{s}}(t-1) - e_{\phi(t)}$
		\State $\mathbf{b}(t) \leftarrow \max\{0,(\mathbf{b}(t-1)-\boldsymbol{\rho})\mathbbm{1}_{t> 1} + e_{\phi(t)} \}$
		\EndFor\\
		\Return $\Phi^{GB}=\{\phi(t)\}_{t=1...T}$
	\end{algorithmic}
\end{varalgorithm}

\subsection{Performance of \ref{alg:modal_balance}}

Recall that the parameter $\gamma$ controls the trade-off between allocation balancing and outcome maximization. Therefore, using historical data the policymaker can tune this parameter to obtain the desired level of employment and allocation balance. In situations in which the payoff/cost of outcomes, wait time, and idle time can all be measured in or converted to a common metric (such as dollars), policymakers might want to set $\gamma$ to the specific value that leads to optimization of that common metric. 

\begin{figure}[h!]
	\begin{center}
		\caption{Trade-off between outcome maximization and allocation balancing for US data (left) and Swiss data (right). The top row shows the results using 2015 data, and the bottom row uses the 2016 data.} \label{fig:tradeoff_curve}
		\includegraphics[width=.4\textwidth]{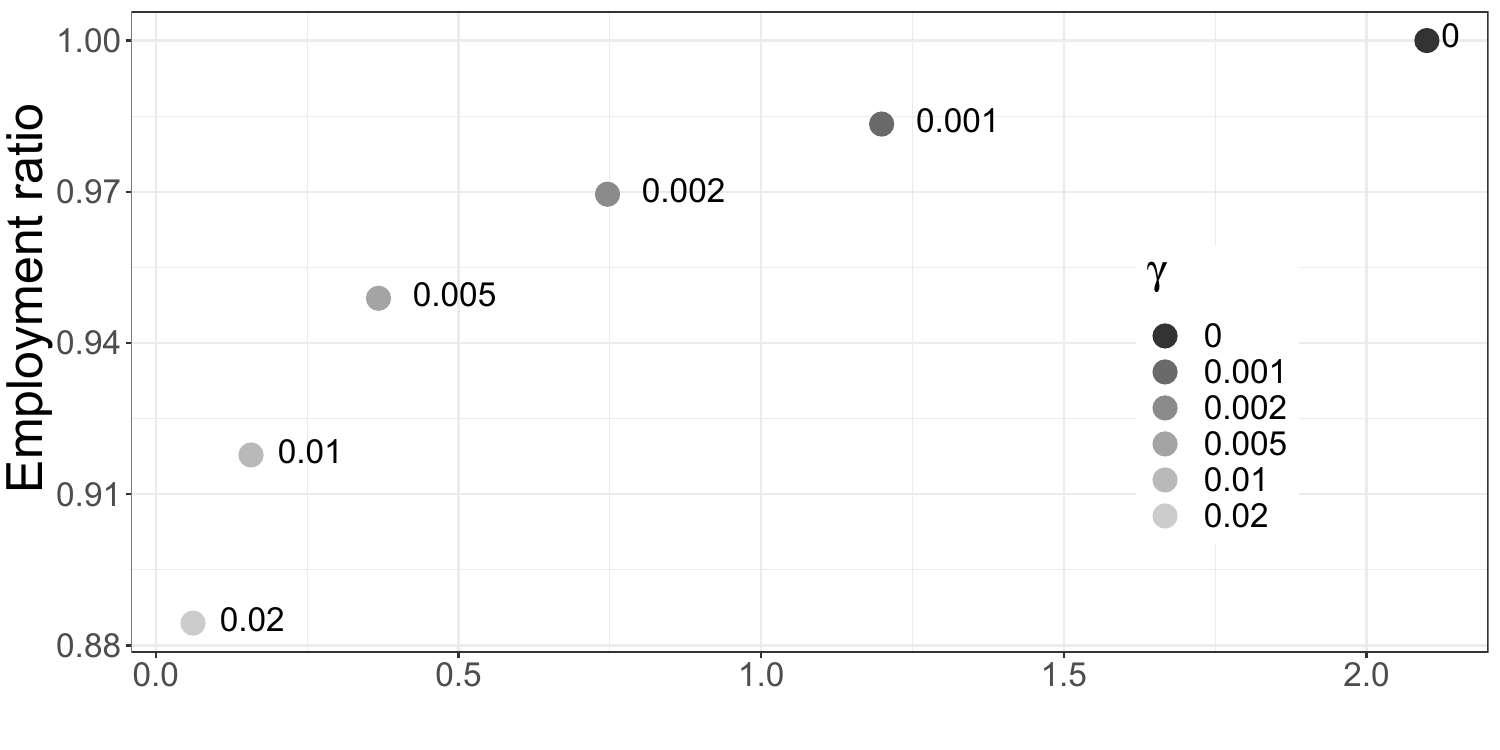}
		\includegraphics[width=.4\textwidth]{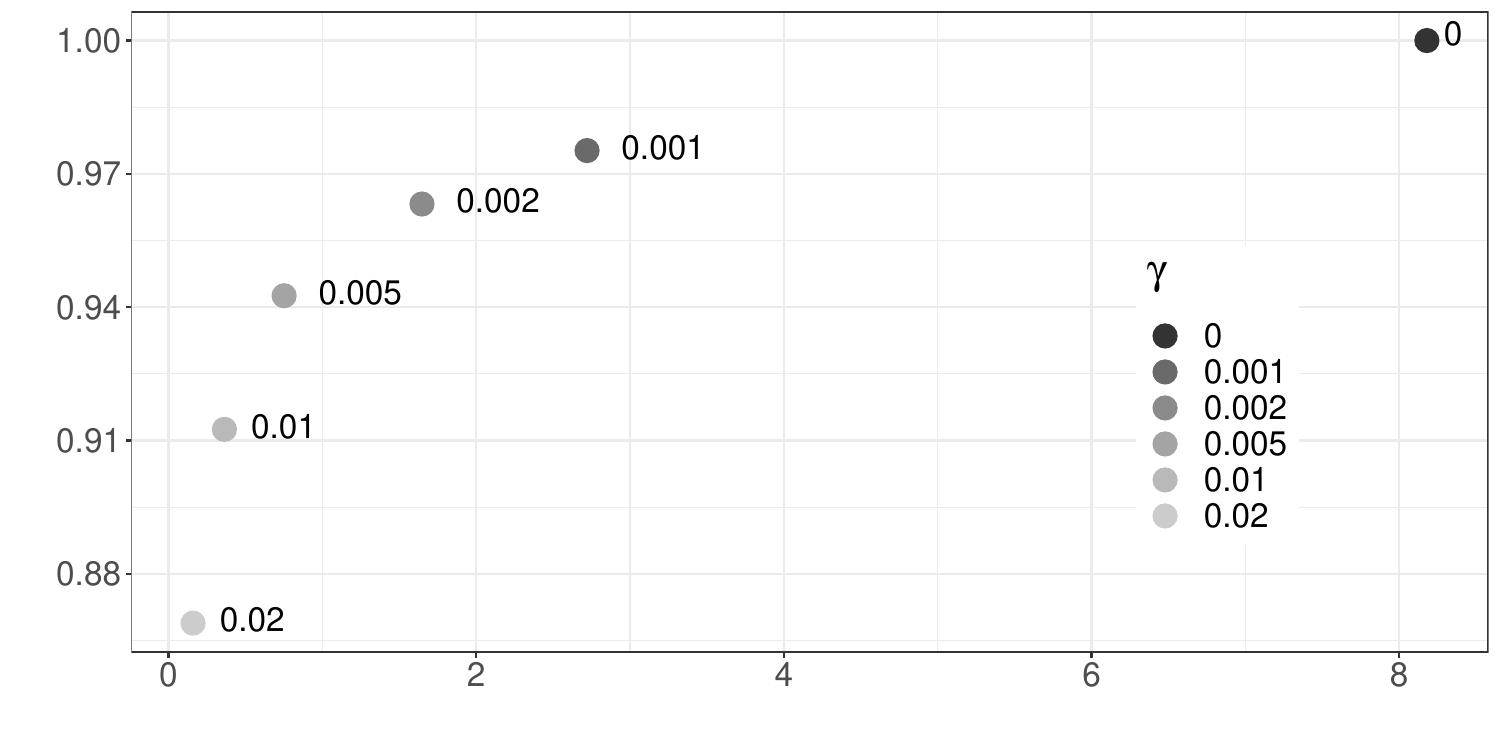}
		\includegraphics[width=.4\textwidth]{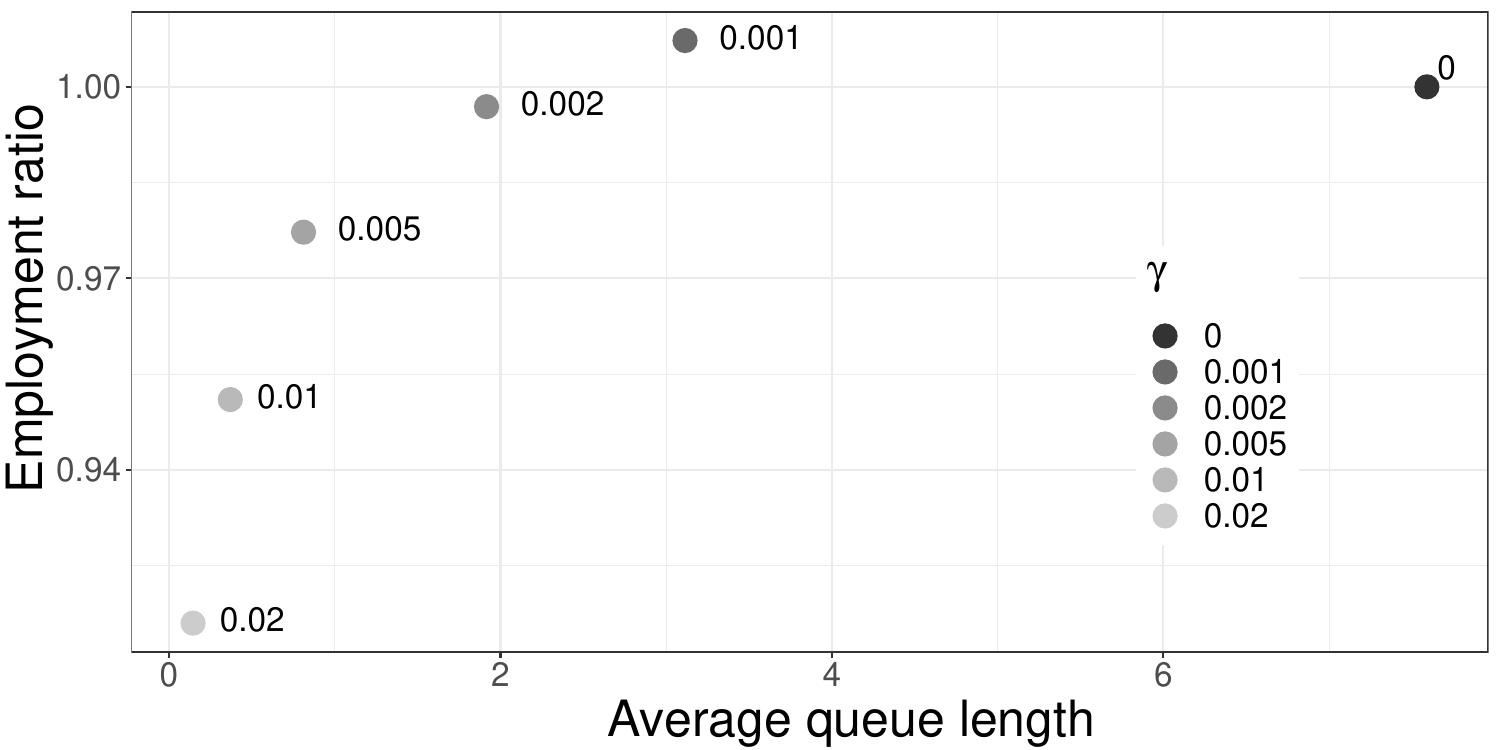}
		\includegraphics[width=.4\textwidth]{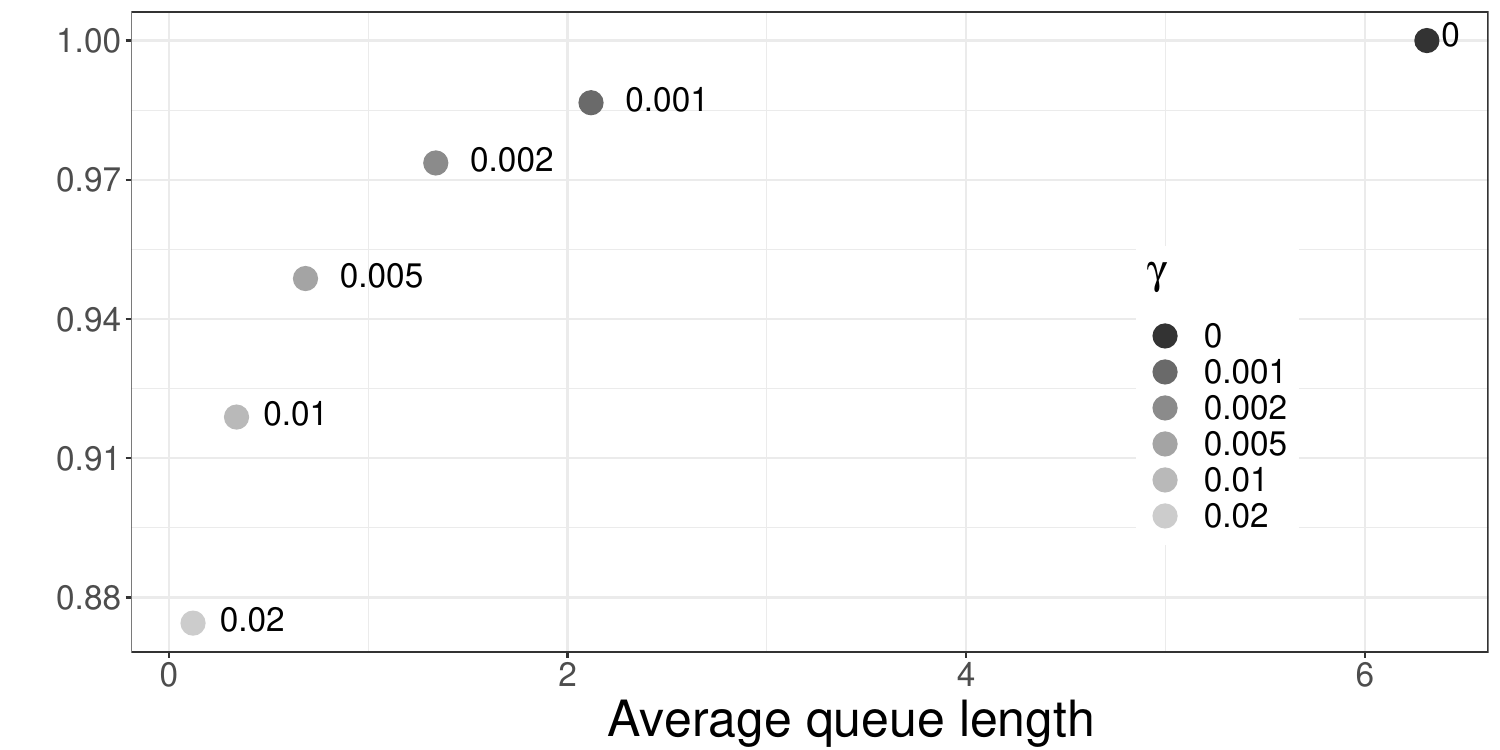}
	\end{center}
\end{figure}

Figure \ref{fig:tradeoff_curve} shows the employment level and average queue length incurred by various values of $\gamma$ on both the 2015 and 2016 arrivals from US and Switzerland. The vertical axis of Figure \ref{fig:tradeoff_curve} shows the employment level under a particular value of $\gamma$ divided by the employment level when $\gamma=0$ (i.e., under pure outcome maximization). The $x$-axis shows the average queue length across affiliates and arrivals.
In practice, we would not know the ``best'' value of $\gamma$ to choose in a given year in hindsight, and would need to base this decision on historical data. Therefore, for the 2016 cohorts, the value of $\gamma$ should be chosen using the top row of Figure \ref{fig:tradeoff_curve} (which uses 2015 data) and the resulting employment and build-up can be seen in the bottom row. Interestingly, on the 2016 US data (Figure \ref{fig:tradeoff_curve} bottom left), the highest employment level is not achieved when $\gamma=0$, but when $\gamma$ is slightly positive, likely due to idiosyncratic non-stationarities in the arrival process.

As can be seen from Figure \ref{fig:tradeoff_curve}, buildup can be dramatically reduced with little loss in employment. The ideal region in Figure \ref{fig:tradeoff_curve} is the top left---where buildup is minimized, and employment is maximized. Based on the top row of Figure \ref{fig:tradeoff_curve}, a policymaker could choose an appropriate value of $\gamma$ to achieve their desired balance of employment versus allocation balancing for the 2016 cohorts. 

\begin{figure}[h!]
	\begin{center}
		\caption{Allocation to nine largest locations over time using \ref{alg:modal_balance} with $\gamma=.005$ for US data (left) and Swiss data (right).} \label{fig:usresults_heuristic}
		\includegraphics[width=.49\textwidth]{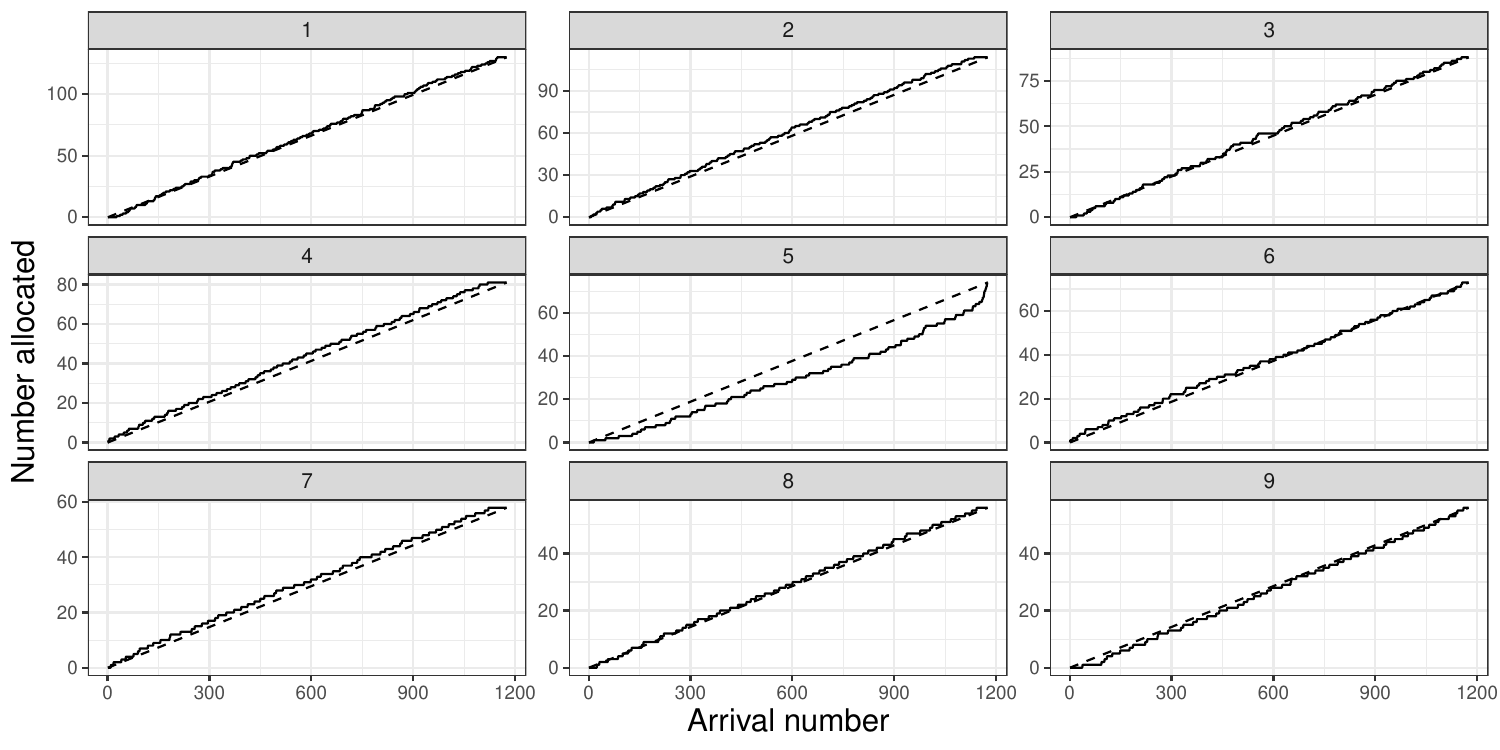}
		\includegraphics[width=.49\textwidth]{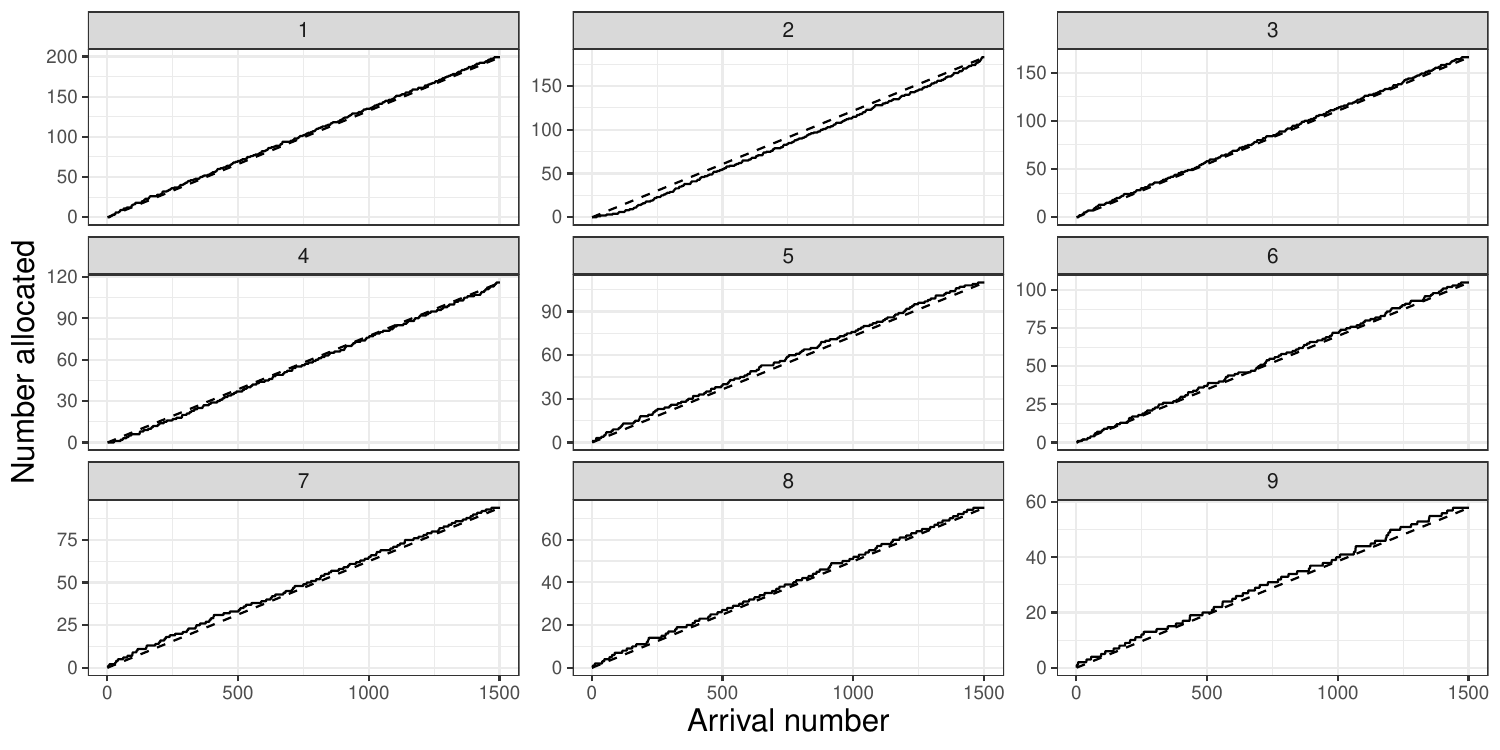}
	\end{center}
\end{figure}

To illustrate these results in greater detail, Figure \ref{fig:usresults_heuristic} shows the allocation of the 2016 arrivals using \ref{alg:modal_balance} with $\gamma=0.005$. This can be compared to Figure \ref{fig:usresults_alg_imbalance}. From visual inspection alone, it is clear that \ref{alg:modal_balance} results in a much more balanced allocation over time. Indeed, the average queue length is less than one.
 Furthermore, the total employment level obtained using \ref{alg:modal_balance} with $\gamma=0.005$ is 98\% of the level obtained with \ref{alg:modal}. On the Swiss data, the results are similar: the average queue length is less tha one
 and the employment level obtained is 95\% of the level obtained under \ref{alg:modal}.

Thus, with little loss in employment, \ref{alg:modal_balance} is able to achieve a highly balanced allocation over time. A balanced allocation results in an even workload for resettlement officers and immediately solves many issues associated with updating capacities over time, discussed in Section \ref{sec:implementation}. Specifically, the capacity of each location could be set according to an upper confidence bound on the number of arrivals, without running the risk of particular locations exceeding their proportionality key by the end of the horizon. 

\subsection{Ancillary Benefit: Increased Exploration}

Although this paper does not focus on the outcome prediction methodology, the prediction and assignment steps are not independent (as discussed in \cite{kasy2023matching}). In this paper, it was assumed that the outcome scores are known. In practice, these outcome scores are estimated from historical data. To use the proposed methods, reliable scores must be determined for every combination of covariates and locations.
If these scores are generated via statistical estimation procedures, maintaining some degree of exploration---assigning similar cases to different locations---is crucial to the resiliency of the estimation procedure given the non-stationarity of the environment. The need for exploration in these situations is a well-known issue and is not unique to the refugee matching context.

A non-stationary constrained contextual bandit framework could be used to formally address this problem. However, without formalizing the bandit version of this problem, we note that \ref{alg:modal_balance} achieves higher levels of exploration than \ref{alg:modal}. Intuitively, because of the balancing component of the objective function in \ref{prob:offline_balance_greedy}, the assignment of a case not only depends on their predicted employment score and the remaining capacity vector, but also depends on the current build-up at each location, effectively adding a degree of randomness to the assignment. 

\begin{figure}[h!]
	\begin{center}
		\caption{Average probability of being assigned to the $k$th ranked location, where locations are ranked at the case-level according to their assignment probabilities.} \label{fig:exploration}
		\includegraphics[width=.7\textwidth]{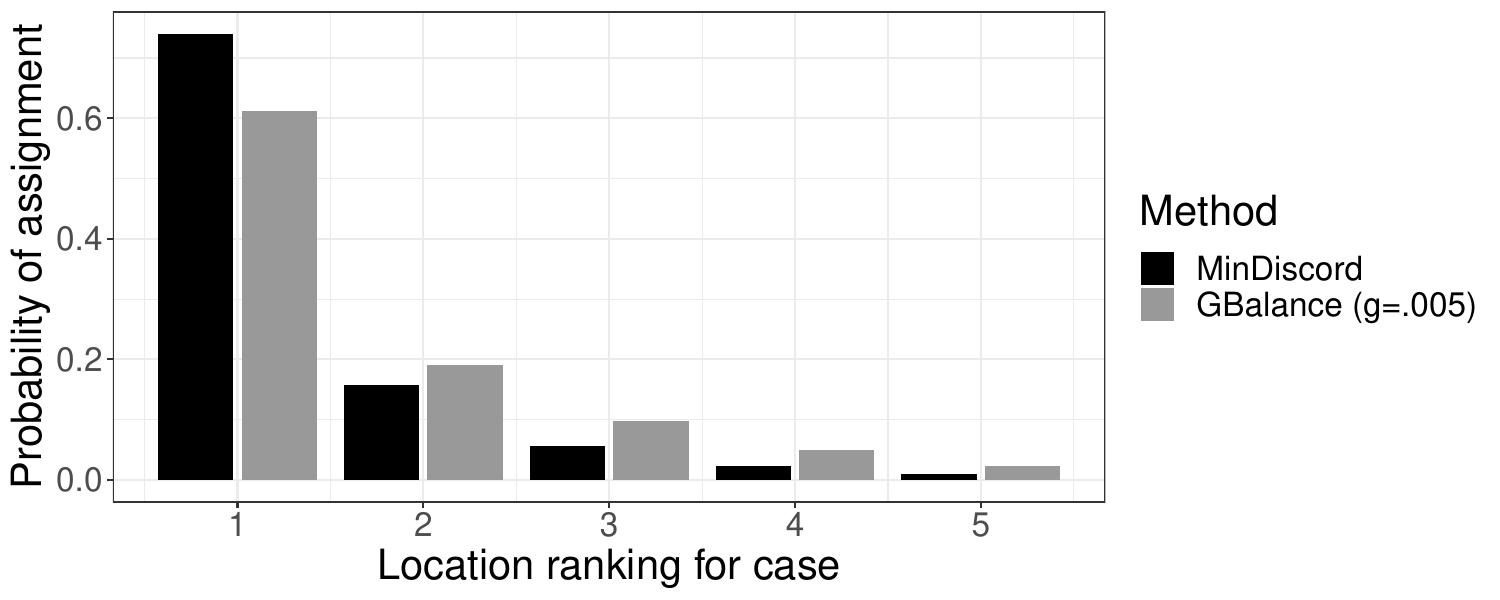}
	\end{center}
\end{figure}

To demonstrate this idea, we run \ref{alg:modal} and \ref{alg:modal_balance} 100 times each for the first 100 arriving US cases in 2016, where the arrival order is randomly permuted in each of the 100 instances. Let case $i$ be the case that arrived $i-$th in the true arrival sequence. In the 100 random instances, they could arrive on any of the 100 days. For each case, we compute the number of times that they are assigned to each location. Let $\ell_{i,1}$ be the location that case $i$ is most often assigned to, $\ell_{i,2}$ be their second most assigned location, etc. Let $n_{\ell_{i,k}}$ be the number of times that case $i$ was assigned to their $k$th most-assigned location out of the 100 instances. Figure \ref{fig:exploration} shows a bar chart of the average value of $n_{\ell_{i,k}}/100$ under \ref{alg:modal} and \ref{alg:modal_balance}. Note that $\sum_k n_{\ell_{i,k}}=100$ for each case $i$. If $n_{\ell_{i,1}}=100$, then $n_{\ell_{i,k}}=0$ for all $k>1$, and case $i$ did not ``explore'' at all. The more uniform the values of $n_{\ell_{i,k}}$, the greater the exploration. 

As shown in Figure \ref{fig:exploration}, under \ref{alg:modal_balance} the average value of $n_{\ell_{i,1}}/100$ is about 0.61 (meaning that a case was assigned to their ``top'' location 61\% of the time), whereas under \ref{alg:modal} the value is about 0.74. Additionally, the average number of unique locations that the same case was assigned to under \ref{alg:modal} was 3.97, versus 5.31 under \ref{alg:modal_balance}. This suggests that \ref{alg:modal_balance} may be preferable to \ref{alg:modal} from a resiliency perspective, naturally maintaining a higher degree of exploration.

\section{Conclusions} \label{sec:conclusion}

This study proposed two assignment algorithms for matching refugees to localities. The first method, \ref{alg:modal}, seeks to maximize the employment scores of all refugees over a horizon by minimizing the probability of disagreement between the online algorithm and an offline benchmark. 
On the Swiss asylum-seeker data used in this study, this method is able to achieve 98\% of the hindsight-optimal employment score. This is a significant improvement over the actual historical assignment, random assignment, and greedy assignment, which achieve 64\%, 65\% and 87\% of the hindsight-optimal employment scores, respectively. Similar results are found using US data. \ref{alg:modal} is currently employed in a multi-year pilot in Switzerland.

However,  \ref{alg:modal}---and any outcome maximizing algorithm---may result in severe periodic imbalance across the localities in the presence of non-stationary arrivals. This creates implementation challenges and an imbalanced workload for the local caseworkers, service providers, and other community members who help each newly arriving family get settled. Furthermore, if local capacities must be revised throughout the year due to larger or smaller arrival numbers than anticipated, imbalance in the allocation over time makes that capacity revision process more challenging. Therefore, we proposed a second assignment algorithm that directly seeks to balance the allocation over time to the localities, while still achieving high employment levels. On the US and Swiss refugee resettlement data used in this study, the allocation balancing method is able to significantly increase balance with little loss in employment.

By all indications, the challenges and scale of forced migration will continue to grow into the future. The methods presented here build upon recent research on outcome-based refugee assignment and could be integrated into refugee resettlement and asylum programs in many host countries---such as the United States, Netherlands, Switzerland, Sweden, and Norway---to help improve the lives of some of the world's most vulnerable populations. 

\ACKNOWLEDGMENT{We acknowledge funding from the Rockefeller Foundation, Schmidt Futures, Charles Koch Foundation, Stanford Institute for Human-Centered Artificial Intelligence, Google.org, and Stanford Impact Labs. The funders had no role in the data collection, analysis, decision to publish, or preparation of the manuscript. We thank Global Refuge and the Swiss State Secretariat for Migration (SEM) for access to data and guidance. The US and Swiss refugee data used in this study were provided under a collaboration research agreement with Global Refuge and the SEM, respectively. These agreements requires that these data not be transferred or disclosed. The authors are Faculty Affiliates of the Immigration Policy Lab at Stanford University and ETH Zurich. This work is part of IPL's GeoMatch project. Finally, we are grateful to Daniel Freund, Jens Hainmueller, and Dominik Hangartner for helpful comments.}

\bibliographystyle{plainnat}
\bibliography{references}

\ECSwitch

\setcounter{endnote}{0}

\ECHead{Online Companion}

\section{Regret of \ref{alg:modal}}\label{app:regret}

The total regret of \ref{alg:modal} is defined as:
$$R(\ref{alg:modal}) :=w(\Phi^{MD})-w(\Phi^*) $$
where $w(\cdot)$ is defined in Equation \eqref{eq:totalscore} and $\Phi^*$ denotes the hindsight-optimal matching (namely, the solution to \ref{prob:offline}$(\mathbf{W})$ where $\mathbf{W}$ is the true test cohort).

Suppose that, at time $t$, an arbitrary online algorithm assigns case $t$ to location $\phi(t)$. The expected \emph{disagreement cost} is defined as
\begin{equation} \label{eq:disagreement_cost}
	\begin{aligned}
		d_t(\phi(t),S_t) = \mathbb{E}[  &\ref{prob:offline}([\mathbf{w}_t;\mathbf{W}(\omega)_{t+1}], \tilde{\mathbf{s}}(t-1))\\
		&-(w_{t\phi(t)} + \ref{prob:offline}(\mathbf{W}(\omega)_{t+1}, \tilde{\mathbf{s}}(t-1)-e_{\phi(t)}))|S_t].
	\end{aligned}
\end{equation}
This is the amount that the online algorithm loses at time $t$ compared to the offline benchmark starting in the same state. 
Lemma \ref{result:regret} below characterizes the expected regret of \ref{prob:offline} in terms of the disagreement probabilities, analogous to Lemma 1 of \cite{vera2020bayesian}.
\begin{lemma}\label{result:regret}
	Let $\phi(t)^{MD}$ be the location that case $t$ is assigned to under \ref{alg:modal}. Under the stage-wise independence and stationarity assumptions,
	$$\mathbb{E}[R(\ref{alg:modal})]  = \sum_{t=1}^T d_t(\phi(t)^{MD},S_t)\leq\delta_{max} \sum_{t=1}^Tq(\phi(t)^{MD},S_t),$$ 
	where $\delta_{max}$ is an upper bound on the disagreement cost.
\end{lemma}
Before proving Lemma \ref{result:regret}, note that in the setting of this paper, the outcome scores are probabilities so the disagreement costs cannot be greater than one. Thus, an upper bound on the expected regret is simply given by $\sum_{t=1}^Tq(\phi(t)^{MD},S_t)$.
As mentioned before, note that \ref{alg:modal} does not minimize the expected disagreement cost at each step, but rather it minimizes the probability of disagreement at each step. Therefore, it minimizes an upper bound on the expected regret (Lemma \ref{result:regret}). Appendix \ref{sec:min-risk} presents an alternative version of the algorithm that does minimize the expected disagreement cost at each step. 

\begin{proof}{Proof of Lemma \ref{result:regret}}
	Let $v(\cdot)$ denote the objective value of a particular optimization problem. Let $\mathbf{W}$ denote the employment scores for a cohort of interest.
	Let $\Psi_{t}^*(\tilde{\mathbf{s}}(t-1)):=\mathbb{E}_{\mathbf{W}_{t}}[v(\ref{prob:offline}(\mathbf{W}_{t}, \tilde{\mathbf{s}}(t-1)))]$ denote the expected sum of outcome scores of the optimal assignment for all cases $t$ onward, starting with capacity vector $\tilde{\mathbf{s}}(t-1)$. 
	
	Given the assignment of case $t$ to $\phi(t)$, the expected \emph{disagreement cost} at time $t$ can be written as:
	$$d_t = \mathbb{E}_{\mathbf{W}_{t+1}}[ \Psi^*_t(\tilde{\mathbf{s}}(t-1))-\left(w_{t \phi(t)}+ \Psi_{t+1}^*(\tilde{\mathbf{s}}(t-1)-e_{\phi(t)})\right) | \mathbf{w}_t].$$
	This is equivalent to Expression \ref{eq:disagreement_cost} in the main text, where we have dropped the arguments in $d_t(\cdot)$ for simplicity.
	
	Note that the expression above relies on the stage-wise independence assumption.
	Rearranged, and showing the same result for case $t+1$, yields:
	$$w_{t \phi(t)} = \mathbb{E}_{\mathbf{W}_{t+1}}[\Psi_t^*(\tilde{\mathbf{s}}(t-1))  | \mathbf{w}_t] - \mathbb{E}_{\mathbf{W}_{t+1}}[\Psi_{t+1}^*(\tilde{\mathbf{s}}(t-1)-e_{\phi(t)}) | \mathbf{w}_t] + d_{t}$$
	$$w_{(t+1) \phi(t+1)} = \mathbb{E}_{\mathbf{W}_{t+2}}[\Psi_{t+1}^*(\tilde{\mathbf{s}}(t)) | \mathbf{w}_{t+1}] - \mathbb{E}_{\mathbf{W}_{t+2}}[\Psi_{t+2}^*(\tilde{\mathbf{s}}(t)-e_{\phi(t+1)})| \mathbf{w}_{t+1}] + d_{t+1}$$
	
	Now note that, by the assumptions of stationarity and stagewise independence: $$\mathbb{E}_{\mathbf{W}_{t+1}}[\Psi_{t+1}^*(\tilde{\mathbf{s}}(t-1)-e_{\phi(t)}) | \mathbf{w}_t] =\mathbb{E}_{\mathbf{W}_{t+1}}[\mathbb{E}_{\mathbf{W}_{t+2}}[\Psi_{t+1}^*(\tilde{\mathbf{s}}(t)) | \mathbf{w}_{t+1},\mathbf{w}_t]].$$
	where $\tilde{\mathbf{s}}(t) = \tilde{\mathbf{s}}(t-1)-e_{\phi(t)}$.
	
	Therefore, the expected sum of the outcome scores for assigning each item via $\Phi$ can be written as:
	$$E_{\mathbf{W}} \left[ \sum_{t=1}^T w_{t \phi(t)} \right] = E_{\mathbf{W}} \left[ \Psi^*_1(\mathbf{s}) \right] + \sum_{t=1}^T d_t,$$
	and thus the expected regret of assignment $\Phi$ is given by
	$$\mathbb{E}[R(\phi)] = \mathbb{E}_{\mathbf{W}} \left[ \sum_{t=1}^T w_{t \phi(t)}-\Psi^*_1(\mathbf{s}) \right] = \sum_{t=1}^T d_t.$$
	This states that the expected regret of an assignment $\Phi$ is equal to the sum of the expected disagreement costs at each time step. Note that this expression for regret applies to any online algorithm $\Phi$.
	
	Therefore, it is natural to develop an online assignment algorithm that minimizes the disagreement cost at each time period. The algorithm \ref{alg:modal} does not minimize the disagreement cost explicitly, but instead minimizes the probability of disagreement. Appendix \ref{sec:min-risk} proposes an online algorithm that minimizes the disagreement cost directly. This algorithm is more computationally taxing than \ref{alg:modal}.
	
	The sum of expected disagreement costs can also be written as 
	$$\mathbb{E}\left[\sum_{t=1}^T d_t\mathbbm{1}\{Q(\phi(t),S_t)\}\right].$$
	where $Q(\phi(t),S_t)$ is the event that assignment $\phi(t)$ disagrees with the offline benchmark at time $t$. If there is no disagreement, then the disagreement cost is zero.
	This expression can be upper bounded by:
	$$\mathbb{E}\left[\sum_{t=1}^T d_t\mathbbm{1}\{Q(\phi(t),S_t)\}\right] \leq \delta_{max} \mathbb{E}\left[\sum_{t=1}^T \mathbbm{1}\{Q(\phi(t),S_t)\}\right] = \delta_{max} \sum_{t=1}^Tq(\phi(t),S_t)$$
	where $\delta_{max}$ is an upper bound on the disagreement costs. In our setting, because the outcome scores correspond to employment probabilities, the  disagreement cost in any time period is upper bounded by one. Thus, the expected total regret can be bounded above by the number of disagreements:
	$$\mathbb{E}[R(\phi)]\leq \mathbb{E}\left[\sum_{t=1}^T\mathbbm{1}\{Q(\phi(t),S_t)\}\right] = \sum_{t=1}^Tq(\phi(t),S_t).$$
	where $q(\phi(t),S_t)$ is again the disagreement probability of action $\phi(t)$ in state $S_t$. \hfill $\square$
\end{proof}

\section{Proof  of Lemma \ref{lemma:equivalence}}\label{app:lemma-proof}

Recall that $b_j(t)= \max\{0,b_j(t-1)-\rho_j\} +z_{tj}$ is the buildup after the assignment at time $t$ but before the processing at time $t$, and $b_j(1)=z_{1j}$. 
In order to prove the lemma, we must show that the following expressions are equivalent:
\begin{equation}\label{eq:balance1}
	\sum_{t=1}^{T}\sum_{j=1}^M \lceil b_{j}(t)-1 \rceil\mathbbm{1}\{b_{j}(t)>1\},
\end{equation}
 and
\begin{equation}\label{eq:balance2}
	\sum_{t=1}^T \sum_{j=1}^M z_{tj}\left\lceil \frac{b_j(t-1)-\rho_j}{\rho_j}\right\rceil\mathbbm{1}\{b_j(t-1)> 0\}
\end{equation}
The first expression is the wait time component of the objective of \ref{prob:offline_balance}. Intuitively, the former expression calculates the wait time cost incurred \emph{per period}, whereas the latter calculates the costs \emph{per case}.
To see this, notice that 
$\lceil b_{j}(t)-1 \rceil$ is the number of cases waiting (not in-process) at time $t$ at location $j$, and $\left\lceil\frac{b_j(t-1)-\rho_j}{\rho_j}\right\rceil$ is the number of periods that case $t$ waits before being serviced. Let $q_{tj}:=\left\lceil\frac{b_j(t-1)-\rho_j}{\rho_j}\right\rceil$.

Now consider the following string of equalities for expression \ref{eq:balance2}:
\begin{align*}
	\sum_{c=1}^T z_{cj}q_{cj}
	=&\sum_{c=1}^T z_{cj}  \sum_{t=c}^{c+q_{cj}-1} 1\\
	=&\sum_{c=1}^T \sum_{t=1}^T z_{cj}  \mathbbm{1}\{t\in \{c,...,c+q_{cj}-1\}\}\\
	=&\sum_{t=1}^T \sum_{c=1}^T z_{cj}  \mathbbm{1}\{t\in \{c,...,c+q_{cj}-1\}\}.
\end{align*} 
Intuitively, the last equality above exactly counts the number of cases that are waiting to be processed at time $t$ at location $j$.
Case $c$ is waiting in time periods $t\in\{c,...,c+q_{cj}-1\}$ and is being processed in time periods $t\in\{c+q_{cj}, c+\left\lceil\frac{b_j(t-1)-\rho_j+1}{\rho_j}\right\rceil-1\}$ (since $\lceil 1/\rho_j\rceil$ is the processing time). Therefore, $\lceil b_j(t) \rceil$---which is the number of cases being processed \emph{or} serviced at time $t$, is equal to $ \sum_{c=1}^T z_{cj}  \mathbbm{1}\{t\in \{c,...,c+q_{cj}-1\}\}$.

Thus, we can write 
\begin{align*}
	&\sum_{c=1}^T z_{cj}  \mathbbm{1}\left\{t\in \left\{c,...,c+q_{cj}-1\right\}\right\}\\
	&= \lceil b_j(t) \rceil-\sum_{c=1}^T z_{cj}\mathbbm{1}\left\{t\in \left\{c+q_{cj},...,c+\left\lceil\frac{b_j(t-1)-\rho_j+1}{\rho_j}\right\rceil-1\right\}\right\}.
\end{align*}
Since only one case can be processed at a time, the term
$$\sum_{c=1}^T z_{cj}\mathbbm{1}\left\{t\in \left\{c+q_{cj},...,c+\left\lceil\frac{b_j(c-1)-\rho_j+1}{\rho_j}\right\rceil-1\right\}\right\}$$ can be at most equal to one for each $t$, and is exactly equal to one unless location $j$ is idle at time $t$. Otherwise, this term is zero. Therefore, 
$$\sum_{c=1}^T z_{cj}\mathbbm{1}\left\{t\in \left\{c+q_{cj},...,c+\left\lceil\frac{b_j(c-1)-\rho_j+1}{\rho_j}\right\rceil-1\right\}\right\}= \mathbbm{1}\{b_j(t)>0\},$$
and thus
$$\sum_{c=1}^T z_{cj}  \mathbbm{1}\left\{t\in \left\{c,...,c+q_{cj}-1\right\}\right\} = \lceil b_j(t) \rceil - \mathbbm{1}\{b_j(t)>0\}=\lceil b_j(t)-1 \rceil \mathbbm{1}\{b_j(t)>1\}.$$
Therefore, the objective function of Problem \ref{prob:offline_balance} can be re-written as 
$$\sum_{t=1}^T\sum_{j=1}^M w_{tj}z_{tj}-\gamma \sum_{t=1}^T \sum_{j=1}^M  z_{tj} \left\lceil\frac{b_j(t-1)-\rho_j}{\rho_j}\right\rceil\mathbbm{1}\{b_{j}(t-1)>0\}.$$ \hfill $\square$

\section{Additional modeling complexities} \label{app:complexities}

This section explores three modeling complexities not present in the Swiss pilot but which may emerge in future implementations: 1) varying case sizes, 2) tied cases, and 3) constraints on free cases.

\subsection{Varying case size} \label{sec:non-unit-case}

The first complexity is varying case sizes, which adds complexity to the algorithm when location capacities are set at the individual level. As noted in the main text, in the Swiss pilot capacities are specified at the case level. Thus, we need not consider varying case sizes (other than noting that $w_{tj}$ is interpreted as the average probability of employment across adults in case $t$). This decision to use case-level capacities is specific to the pilot implementation. In future Swiss implementations and implementations in other contexts, it may become necessary for the algorithm to take into account individual-level capacities. As noted in \cite{ahani2021dynamic}, US resettlement agencies are given target individual-level quotas for each affiliate by the US State Department. The goal of the agencies is to resettle cases so that the total number of individuals remains within the quotas, although exceeding the quotas by less than 10\% is acceptable. 
Thus, evaluating the potential extensions of the proposed algorithm, and its performance under circumstances where cases vary in size and capacities are established at the individual level, is of interest.

\ref{alg:modal} can be run exactly as before with a slightly modified version of \ref{prob:offline}. In particular, we now solve the following offline problem at time $t$:

\OneAndAHalfSpacedXI
\begin{equation} \label{prob:offlineVS} \tag{\textsc{OutcomeMaxVS}}
	\begin{aligned}
		\maxA_{\mathbf{Z}} &\sum_{t'=t}^{T}\sum_{j=1}^M w_{t'j} z_{t'j}\\
		s.t. 	& \sum_{j=1}^M z_{t'j} = 1 \;\;\; \forall \; t' \in \{t,...,T\}\\
		& \sum_{t'=t}^T z_{t'j} l_{t'} \leq  \tilde{s}_j(t-1) \;\;\; \forall \; j \in [M]\\
		& \mathbf{Z}\in\{0,1\}^{(T-t+1) \times M} \\
		&
	\end{aligned}
\end{equation}
\DoubleSpacedXI
where $l_t$ is the size of case $t$, and $\tilde{s}_j(t-1)$ is again the remaining capacity at location $j$ at the start of time period $t$, however this capacity now applies at the individual level. In the experiments in the main text, the capacities were set so that $s_j$ is equal to the total number of cases that were assigned to location $j$ in the data. When dealing with cases of varying size, we allow for a realistic amount of slack in the capacities. Otherwise, if $s_j$ were set to be exactly equal to the number of individuals assigned to location $j$ historically, \ref{prob:offlineVS} would not be feasible for many randomly sampled arrival trajectories. 
Therefore, when we consider varying case sizes, we define $s_j$ to equal 110\% of the number of individuals actually assigned to location $j$.
This modification is grounded in reality as US resettlement agencies allow for assignments that use up to 110\% of the planned capacity at each affiliate. In implementing \ref{alg:modal}, we may still encounter sample trajectories that are infeasible. To deal with this, any future arrival trajectories that are sampled at time $t$ that are not feasible are discarded. Therefore, the parameter $K$ indicates the number of \emph{feasible} sample paths that will be sampled for each assignment.

\begin{figure}
	\centering
	\includegraphics[width=.6\textwidth]{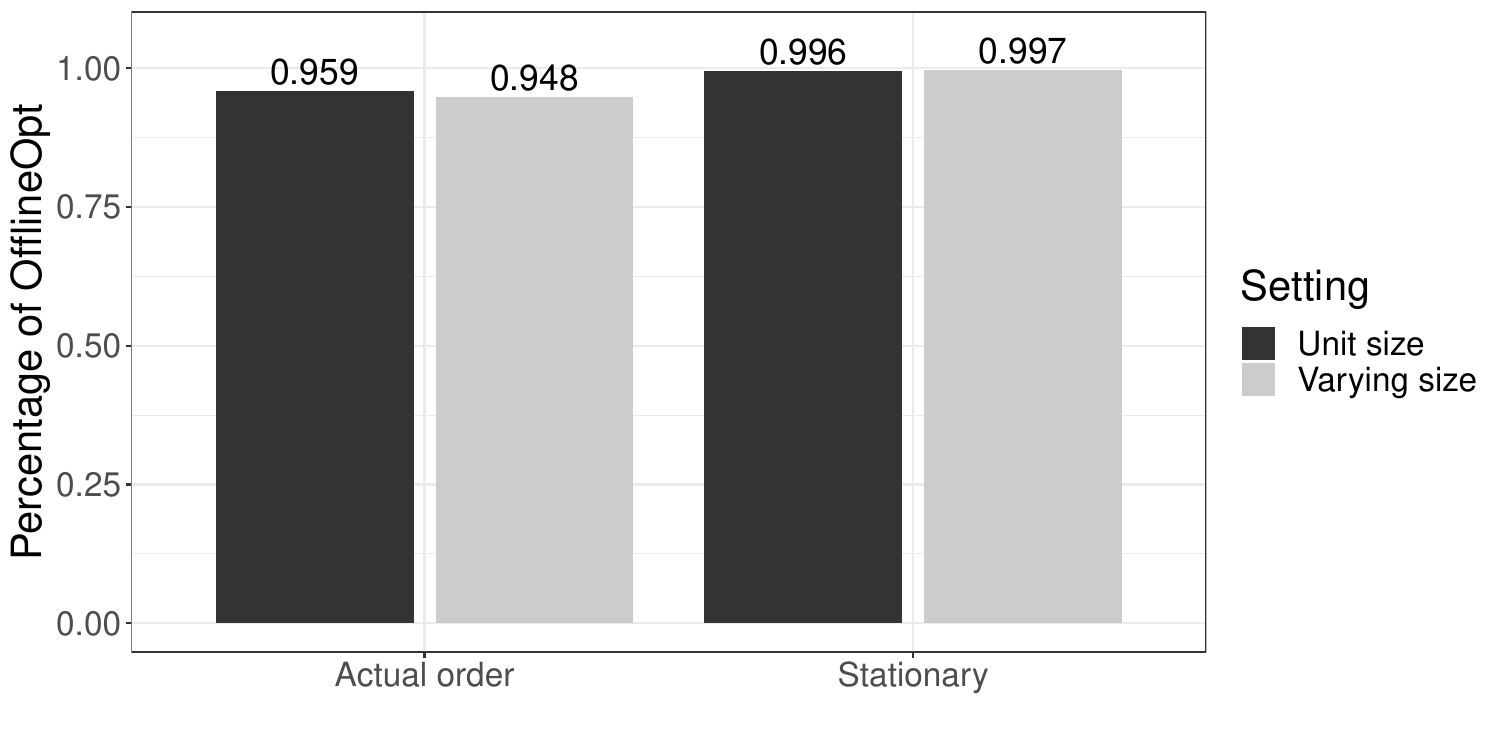}
	\caption{Performance, relative to \textsc{OfflineOpt}, of \ref{alg:modal} using unit case sizes and varying case sizes on the actual 2016 arrivals (left) and randomly drawn arrivals that mimic a stationary process (right). For the varying case size setting, both \ref{alg:modal} and \textsc{OfflineOpt} are allowed to use up to 110\% of the capacity at each affiliate.} \label{fig:varying-case-size}
\end{figure}

Figure \ref{fig:varying-case-size} compares the performance of \ref{alg:modal} when case size is assumed to be equal to one (as in the main text) versus varying case size for the 2016 US data. The offline optimal solution, \textsc{OfflineOpt}, refers to the hindsight-optimal solution obtained by solving \textsc{OutcomeMax} (or \textsc{OutcomeMaxVS}) for the entire horizon in the unit (or varying) case size setting.
The left-hand side of Figure \ref{fig:varying-case-size} shows the performance on the real 2016 data where cases arrive in their actual order. 
The right-hand side considers a setting in which the data is constructed to mimic a stationary process. As described in the main text, to mimic a stationary process the 2016 cohort's arrival sequence is randomly permuted, and $\mathcal{A}$ is set to be the 2016 arrival cohort. Therefore, both the arrival and sampling distributions are identical, and the arrival sequence is stationary. With stationary data, the optimality percentage remains close to 100\% even with the inclusion of varying case sizes. In the context of the actual 2016 data, the optimality percentage decreases from 95.9\% to 94.8\%. Hence, while the inclusion of varying case sizes does induce some performance degradation relative to the offline optimal solution, this loss is mimimal.

\subsection{Family ties}  \label{sec:tied-cases}

As noted in the main text, only free cases (i.e., cases that had no preexisting family ties) were included in the empirical demonstrations of our methods. This reflects the current pilot in Switzerland, where the scope of the pilot is limited to free cases, and separate capacity was reserved for the pilot group so that there is no interference between the non-pilot and pilot cases. Thus, the model and algorithm presented in Sections \ref{sec:problem} and \ref{sec:outcome_max} most closely reflect the pilot implementation.

That being said, future implementations of the algorithm (in Switzerland or in other countries) could include both free and tied cases. In reality, the extent to which free and tied cases share capacity can vary across settings. At one extreme is a situation in which there is no advance information on how many arrivals will have family ties (and where those family ties will be), such that assignment of both free cases and family-tie cases must be made according to a shared capacity system, where they are ``competing'' for the same slots. At the opposite end of the spectrum, arrivals with family ties can be predicted, determined, or selected beforehand. This would allow the location capacities to be separated out for free cases (i.e. the total quotas net of those dedicated to the family-tie cases), as in how the empirical demonstrations have been presented in this paper. In the United States and other countries, the reality will be somewhere in between these two extremes, as it is generally possible to use recent trends and knowledge of prior arrivals who have indicated they have ``trailing'' family members to control or project (with a degree of uncertainty) the number of arrivals who will have family ties and where those ties are.

How tied cases should factor into buildup may also vary by context. At one extreme, those cases would contribute to build-up in the exact same manner as free cases, which would then have implications for balancing the dynamic assignment decisions made for free cases. At the other extreme, the family networks and relationships that family-tie cases have at their receiving locations would eliminate the processing time and costs imposed on resettlement resources and hence result in negligible or no impact on build-up. The reality is somewhere in the middle.

Several approaches could enable the direct integration or absorption of the assignment of family-tie cases into the algorithms presented in this study. The most straightforward way that we explore in this section, which does not use any information about the future arrival of tied cases, is to simply include the family tie cases in the algorithm while forcing their assignment to their predetermined locations, assuming that quotas are shared and all cases contribute to build-up in the same manner. In this straightforward extension, the historical data from which the algorithm samples random arrival trajectories should include both free and tied cases. Intuitively, this approach will likely work well when the number of tied cases that are predetermined for each location is relatively constant year-to-year (i.e., when tied case arrivals are stationary). In what follows, we demonstrate the performance of this simple extension. We leave as an area for future research the development of more robust algorithms that could take distributional knowledge or forecasts of tied cases into account.

\begin{figure}
	\centering
	\includegraphics[width=\textwidth]{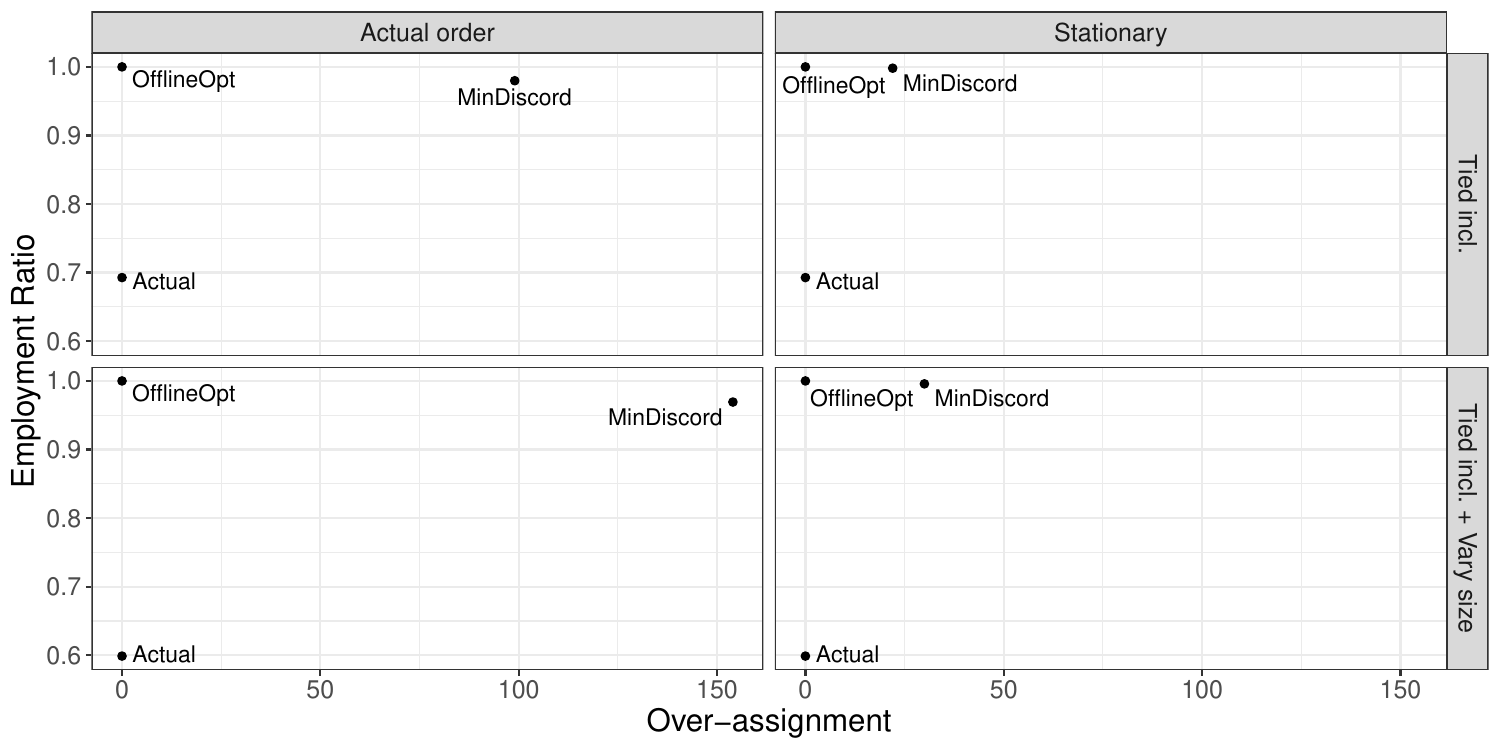}
	\caption{Performance of \ref{alg:modal} with the inclusion of cases with family ties (left) and family ties plus varying case sizes (right). The $y$-axis shows the average employment achieved as a percentage of the offline optimal solution, and the $x$-axis shows the number of individuals that were assigned in excess of capacity.} \label{fig:including-varying-size-tied-cases}
\end{figure}

Note that when tied cases are included, there is always a chance of over-allocation, as tied cases can arrive at any point in the horizon and must be allocated to their predestined location, even if the location is out of capacity. Thus, comparing the performance of algorithms requires comparing the average employment rate and degree of over-allocation (and potentially also the degree of imbalance). Figure \ref{fig:including-varying-size-tied-cases}  demonstrates the performance of \ref{alg:modal} when tied cases (with unit and varying case sizes) are included for the US data. Overall in this expanded dataset, 70\% of cases had pre-existing family ties in 2016 and 75\% in 2015. 

Figure \ref{fig:including-varying-size-tied-cases} demonstrates the performance of \ref{alg:modal} with tied cases included (top row), and when we consider varying case sizes in addition to the inclusion of tied cases (bottom row). In Figure \ref{fig:including-varying-size-tied-cases}, over-assignment is shown along the $x$-axis of each panel and the average employment relative to the offline optimal solution is shown along the $y$-axis. Over-assignment is measured as the total number of individuals that were assigned to locations in excess of the locations' capacity. Recall that when varying case sizes are considered, we assume that the capacity of a location is 110\% of the number of individuals actually assigned to that location in the data.

Figure \ref{fig:including-varying-size-tied-cases} shows the performance on the real-world 2016 arrival data (left column) as well as the performance on synthetically created (as described above) stationary arrival data (right column). The purpose of this is to better understand what portion of the performance difference beween \ref{alg:modal} and \textsc{OfflineOpt} is due to nonstationarity in the data versus the online nature of the problem. 

In this implementation, \ref{alg:modal} does not use any knowledge of the arrival of tied cases, contributing to significant over-allocation. In practice, the resettlement agency is likely to have some knowledge regarding the arrival of future tied cases that could be utilized to limit the resulting degree of over-allocation. Future research could explore methods expressly aimed at minimizing over-allocation.

Lastly, with the inclusion of tied cases, we observe a similar trade-off between allocation balance and employment, as displayed in Table \ref{tab:varying-size-tied-cases-balance}. Notice that increasing $\gamma$ in \ref{alg:modal_balance} not only decreases average wait time as it was designed to do, but also decreases over-assignment. As in the main text, increasing $\gamma$ also has an impact on average employment. This loss in efficiency is larger when tied cases are present (e.g., the average employment when $\gamma=0.005$ is 91\% of the average employment when $\gamma=0$, and 88\% of the average employment obtained by \textsc{OfflineOpt}). Intuitively this makes sense: in the presence of tied cases, achieving a balanced allocation necessitates assigning a greater proportion of free cases to locations where fewer tied cases have arrived. Thus, for example, if tied cases tend to be predestined to urban centers where employment rates are generally higher, this means that more free cases will be assigned to the non-urban centers (which also explains the reduction in over-assignment).

As a reference, when incorporating tied cases and varying family sizes, the actual historical allocation yields an average queue length of 24.5. Thus, when $\gamma\geq 0.001$, \ref{alg:modal_balance} achieves an average queue length that is less than the status quo.

\begin{table}
	\centering
	\begin{tabular}{cccc}
	$\gamma$ & Average employment ratio& Over-assignment & Average wait time \\
	\hline
	0 & 0.97 & 154 & 38.0\\
	0.001 & 0.94 & 28 & 23.0 \\
	0.005 & 0.88 & 23 & 18.1
	\end{tabular}
\caption{This table shows 1) Average employment, as a percentage of \textsc{OfflineOpt}, 2) Over-allocation, computed as the total number of arrivals that are assigned to locations beyond the location's capacity, and 3) the average wait time achieved by \ref{alg:modal} ($\gamma = 0$) and \ref{alg:modal_balance} with $\gamma \in \{0.001,0.005\}$ using 2016 US data including tied cases and varying family sizes.} \label{tab:varying-size-tied-cases-balance}
\end{table}

\subsection{Free case constraints}  \label{app:free-case-constraints}

As noted in the main text, the pilot implementation in Switzerland is limited in scope to free cases that can be assigned to any canton. Thus, the model and algorithm presented in Sections \ref{sec:problem} and \ref{sec:outcome_max} most closely reflect the pilot implementation. In reality, some free cases have idiosyncracies such as educational or medical considerations that may prohibit their assignment to certain locations. 
Practically, incorporating such constraints into the assignment problem is straightforward: if case $t$ cannot be assigned to location $j$, we can set $w_{tj}=-M$ for some arbitrary constant $M>1$ and proceed with \ref{alg:modal} or \ref{alg:modal_balance} as described in the main text.

We do not have access to reliable historical data on these idiosyncratic constraints and thus are unable to evaluate how they might impact the algorithms' performance. However, one extreme modeling assumption would be to treat these cases as tied cases and assume that there was only one location to which they could have been assigned. This would effectively increase the proportion of tied cases present in the data, and we can expect similar but more extreme degredation in the algorithms' performance as discussed in the above section.

\section{Batching} \label{sec:batching}

In the dynamic formulation of the assignment problem presented above, each item is observed and must be assigned one by one. This reflects the assignment procedure in countries such as Switzerland and the Netherlands. However, in other contexts, there may be periodic (e.g. weekly, monthly) cohorts of refugee arrivals that can be assigned in batches rather than on a purely one-by-one basis. This reflects the assignment procedure in the US where assignments are made in weekly batches. This section extends the dynamic assignment mechanisms presented above to a batching context.

First, we note that batching cannot hurt the performance of an online algorithm, since the batches can simply be ignored and the units within batches could be assigned one-by-one using the proposed online assignment algorithms. Batching, therefore, only presents an opportunity for improved efficiency. Using the 2016 US data as an example, we saw that \ref{alg:modal} achieves 95\% of OfflineOpt. Thus, there is an efficiency gap that batching could help to close. We also note that when arrivals are stationary, the optimality gap is only 0.5\%, so most of the efficiency gap is caused by the non-stationarity of arrivals. Because batching allows the algorithm to ``see into the future'' more than a purely online method, it is possible that batching could help alleviate some of the inefficiencies caused by non-stationarity. However, we note that weekly batching still does not allow for too much insight into the future.

We propose two methods for applying our proposed algorithms to a batched setting. 
The first method assigns each case in the batch simultaneously, resulting in the largest efficiency gains. After fixing the arrivals in the current batch,
we again sample $K$ random trajectories of arrivals for the remaining horizon
$$\{\mathbf{z}^*_l\}_{l\in \{t,...,t+B_t\}} = \text{mode} \left(\{\mathbf{z}_{l}^{1}\}_{l\in\{t,...,t+B_t\}},...,\{\mathbf{z}_{l}^{K}\}_{l\in\{t,...,t+B_t\}}\right)$$
In the pure online setting without batching, this reduces to Method \ref{heur:modal-balancing}. In the pure online setting, each vector $\mathbf{z}_t^{k}$ is effectively one dimensional (since it only contains one positive element, it can be mapped to the one-dimensional domain $\{1,...,M\}$). Therefore, the disagreement probabilities could be estimated with a reasonable number of samples. However, with batching, the assignment decision is comprised of $B_t$ assignment vectors, and thus can be mapped to a $B_t$-dimensional domain ($\{1,...,M\}^{B_t}$). Thus, given the combinatorial complexity of the assignment decision with batching, obtaining reasonable estimates of the disagreement probabilities potentially requires many more sample trajectories. When $B_t$ is small, this method is likely to be tractable.

When $B_t$ is large, we propose either assigning cases within the batch one-by-one, or by breaking down the batch into smaller batches where the method above can be applied. Notice that even if cases within a batch are assigned one-by-one, knowledge of the employment scores for the entire batch is nonetheless helpful, and improves the performance of \ref{alg:modal} and \ref{alg:modal_balance}. Instead of randomly sampling the entire future horizon, we can fix the employment scores for the remaining cases within the batch, and therefore only need to randomly sample the horizon \emph{after} the last case in the batch.

\section{Run-time of Offline Problems} \label{app:scale}

This section presents the run-time performance of the proposed approaches. Although we can solve \ref{prob:offline} and \ref{prob:offline_balance_greedy} efficiently as linear programs, run-time considerations may still exist if the horizon is large. Each of the $K$ for-loops in \ref{alg:modal} or \ref{alg:modal_balance} could operate in parallel. Thus, our focus specifically lies on the time required to solve a single instance of the offline problems \ref{prob:offline} (equivalent to the run-time of \ref{prob:offline_balance_greedy}) and \ref{prob:offline_balance}.

For these run-time results, we set the number of locations equal to 30 and randomly create outcome matrices of size $T\times 30$ where each element is independently drawn from a uniform distribution on $[0,1]$. Capacities for the locations are set randomly such that their sum is equal to $T$. Furthermore, to obtain the best possible performance of \ref{prob:offline_balance}, we relax the nonlinearity in the objective function and linearize the Max$(\cdot)$ operator in the constraints. Specifically, we test the run-time of the following modified version of \ref{prob:offline_balance}:

\OneAndAHalfSpacedXI
\begin{equation} \label{eq:fullbalance-app}
	\begin{aligned}
		\maxA_{\mathbf{Z},\mathbf{b}, \tilde{\mathbf{b}}} &\sum_{t=1}^{T}\sum_{j=1}^M w_{tj}z_{tj} - \gamma  \sum_{t=1}^{T}\sum_{j=1}^Mb_{j}(t)\\
		\text{s.t.} 	& \sum_{j\in[M]} z_{tj} = 1 \;\;\; \forall \;\;\; t \in [T]\\
		& \sum_{t\in[T]} z_{tj} = s_j\;\;\; \forall \;\;\; j \in [M]\\
		& \tilde{b}_{j}(t) \geq b_{j}(t-1)-\rho_j \;\;\; \forall \;\;\; t \in [T],\; j\in[M]\\
		& \tilde{b}_j(t) \geq 0 \;\;\; \forall \;\;\; t \in [T],\; j\in[M]\\
		& b_j(t) = \tilde{b}_j(t)+z_{tj}  \;\;\; \forall \;\;\; t \in [T],\; j\in[M]\\
		& \mathbf{Z}\in\{0,1\}^{T\times M}
	\end{aligned}
\end{equation}
\DoubleSpacedXI

In Problem \ref{eq:fullbalance-app}, the auxiliary variables $\tilde{b}_j(t)$ will always be equal to $\max\{0, b_j(t-1)-\rho_j\}$ in the optimal solution. Hence, we can write $b_j(t)$ as $b_j(t)=\tilde{b}_j(t)+z_{tj}$.

\begin{figure}
	\centering
	\includegraphics[width=.6\textwidth]{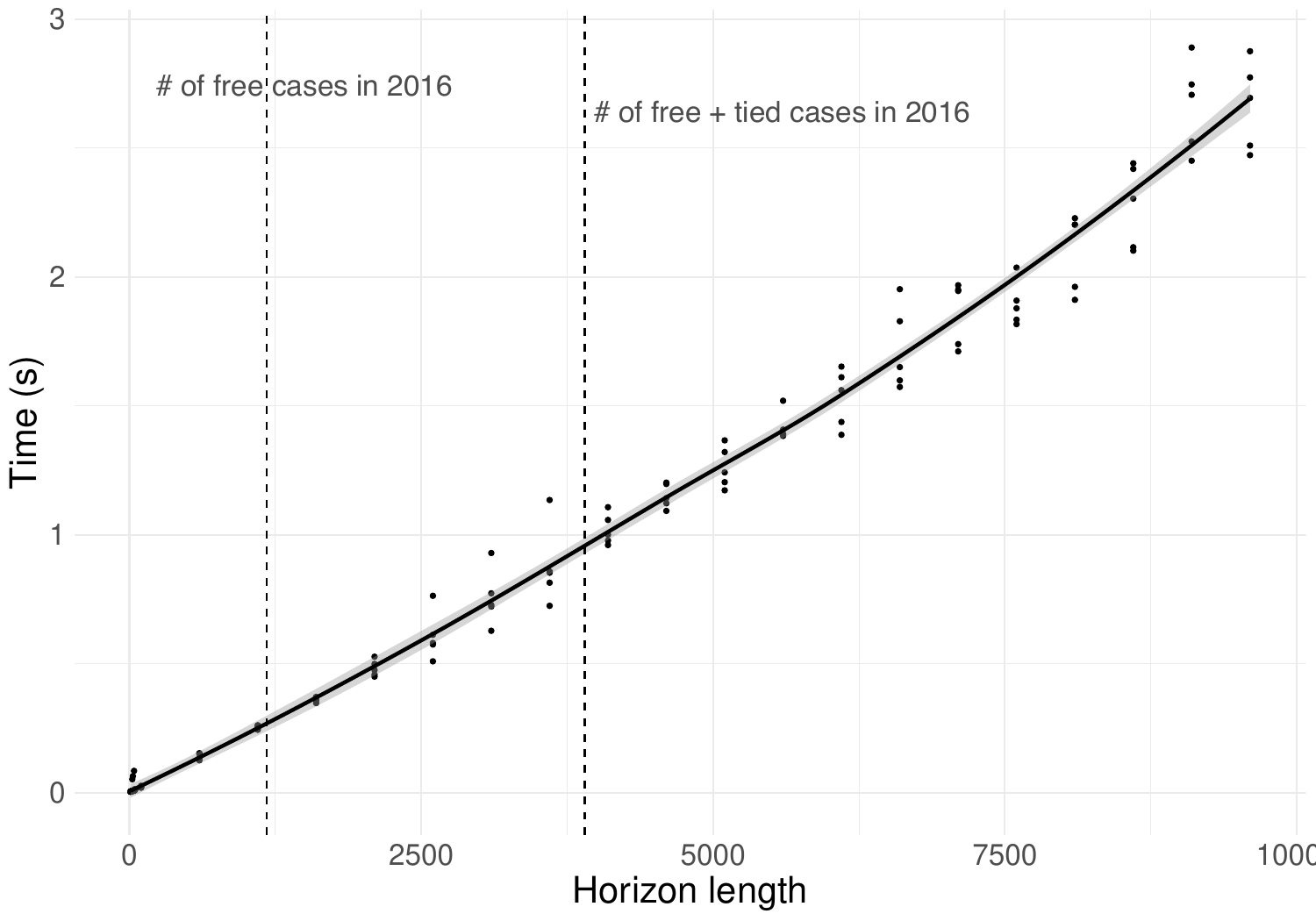}
	\caption{Run-time of \ref{prob:offline} (and \ref{prob:offline_balance_greedy}) as the horizon length increases. For each horizon length, we solve five random instances, shown by the dots on the graph. The blue line illustrates the curve obtained by LOESS smoothing. Run-times are obtained using Gurobi on an iMac with an Apple M1 chip and 16GB of RAM.} \label{fig:run-time}
\end{figure}

Table \ref{tab:balance-runtime} shows that even with very short horizons (of length up to 45), solving Problem \ref{eq:fullbalance-app} is time-intensive, making this approach prohibitively slow for realistic problem sizes that have horizon lengths in the thousands. 
The large jump in run-time when going from a horizon of length 30 to 35 in Table \ref{tab:balance-runtime} can be explained by the fact that there are 30 locations. Namely, when there are fewer people than locations, most locations have a capacity of zero or one, resulting in a relatively simple solution. When this no longer holds, the run-time rapidly increases.
By comparison, Figure \ref{fig:run-time} shows that \ref{prob:offline} (and \ref{prob:offline_balance_greedy}) is fast even as the horizon grows to $T=10,000$.

\begin{table}[h!]
	\centering
	\begin{tabular}{l|cccccccc}
	Horizon length & 10 & 15 & 20 & 25 & 30 & 35 & 40 & 45 \\
	Time (s) & 0.013 & 0.024 & 0.035 & 1.18 & 1.65 & 9.23 & 24.7 & 123 
	\end{tabular}
\caption{Run-time of \ref{prob:offline_balance} for horizons of length 10-45, averaged across 5 random instances. Run-times are obtained using Gurobi on an iMac with an Apple M1 chip and 16GB of RAM.} \label{tab:balance-runtime}
\end{table}

\section{Additional Computations} \label{sec:additional_computations}

\begin{figure}[ht!]
	\begin{center}
 \label{fig:usresults_imbalance_random}
		\includegraphics[width=.47\textwidth]{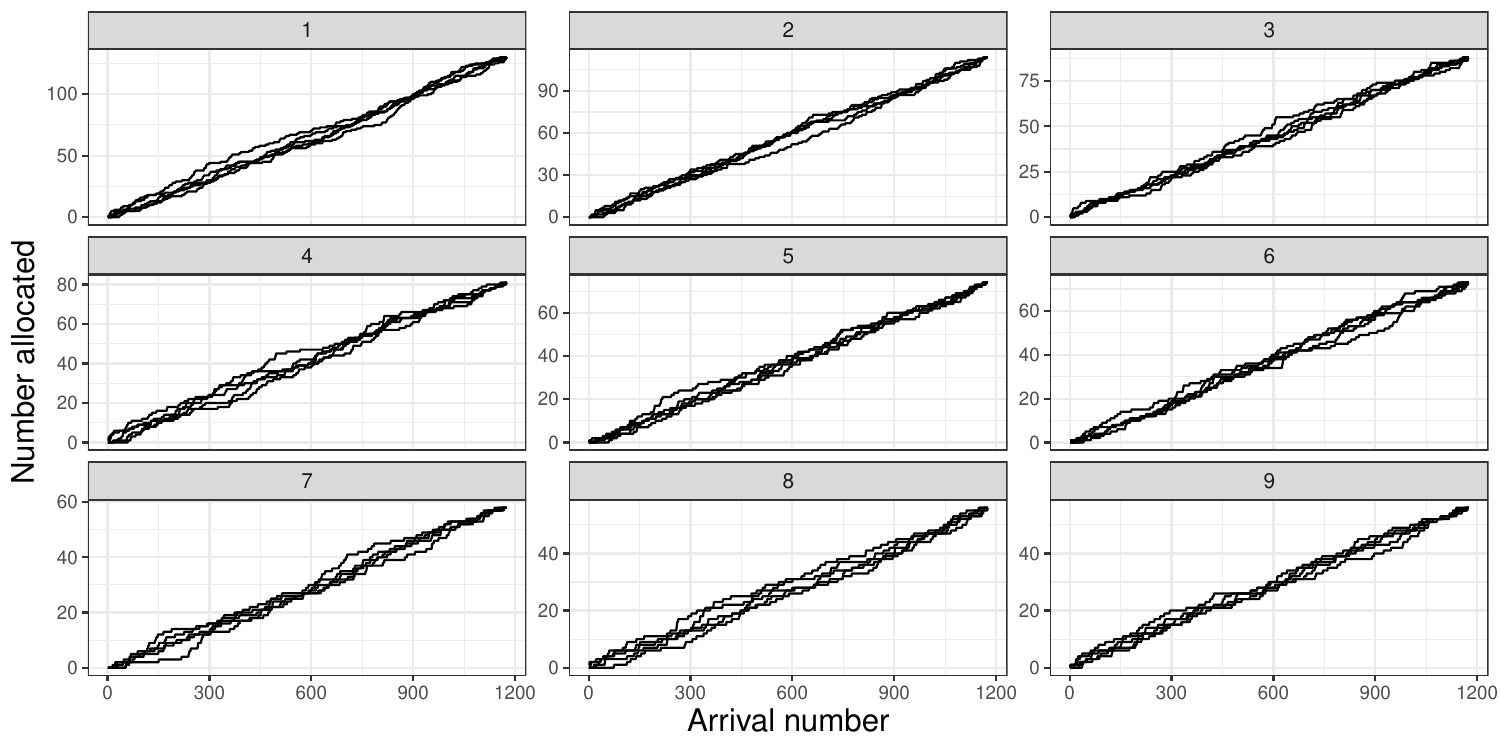}
				\caption{Allocation to top 9 locations over time using \ref{alg:modal} when the arrival and sampling distribution are IID (both drawn at random from the 2016 test cohort). Results are shown across five instances.}
	\end{center}
\end{figure}

Figure \ref{fig:usresults_imbalance_random} shows the allocation over time of \ref{alg:modal} for five random instances where the arrival order of the US test cohort ($N=1,175$) is randomly permuted and the sampling pool is taken to be the 2016 arrivals. Thus, the arrival and sampling distribution are IID. Although Figure \ref{fig:usresults_imbalance_random} is much more balanced than Figure \ref{fig:usresults_alg_imbalance}, random deviations from a balanced allocation can still occur. Across the five instances shown, the average queue length is 1.7. This is similar to the average queue length achieved under a random allocation rule. Thus, even in these cases, explicitly enforcing balance is still beneficial. 

\begin{figure}[ht!]
	\centering
	\includegraphics[width=.45\textwidth]{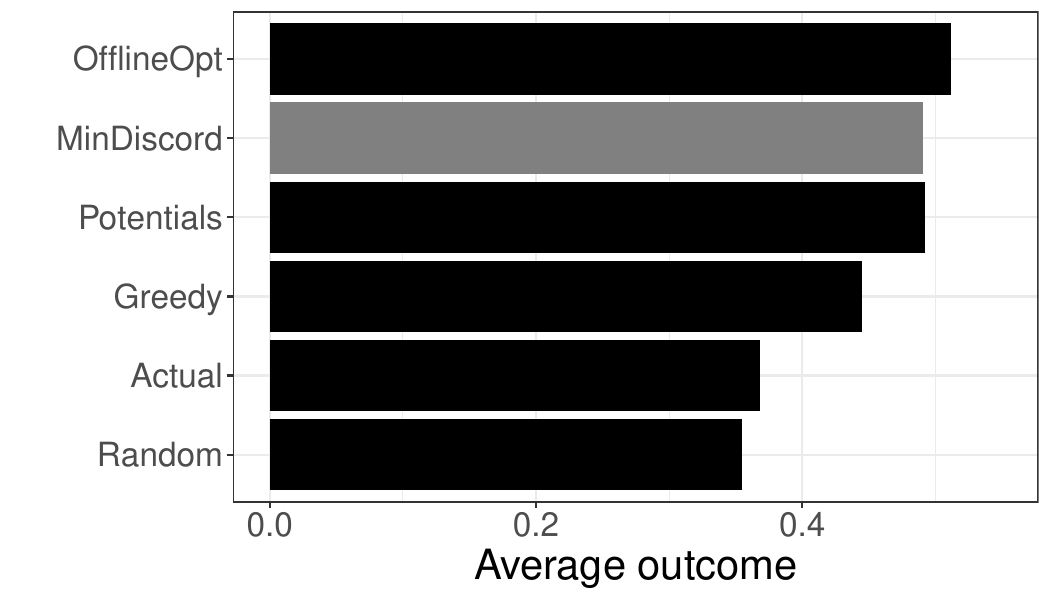}
	\includegraphics[width=.45\textwidth]{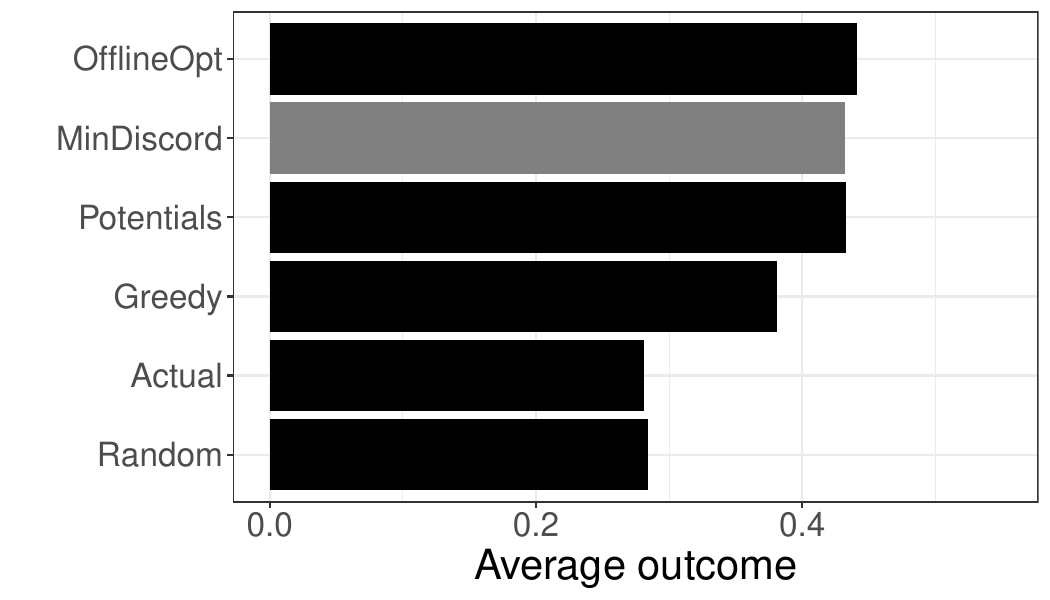}
	\caption{Performance of all outcome maximizing algorithms, including the ``potentials'' method introduced by \cite{ahani2021dynamic} on US data (left) and Swiss data (right).}\label{fig:potentials_comparison}
\end{figure}

\begin{figure}[ht!]
	\centering
	\includegraphics[width=.45\textwidth]{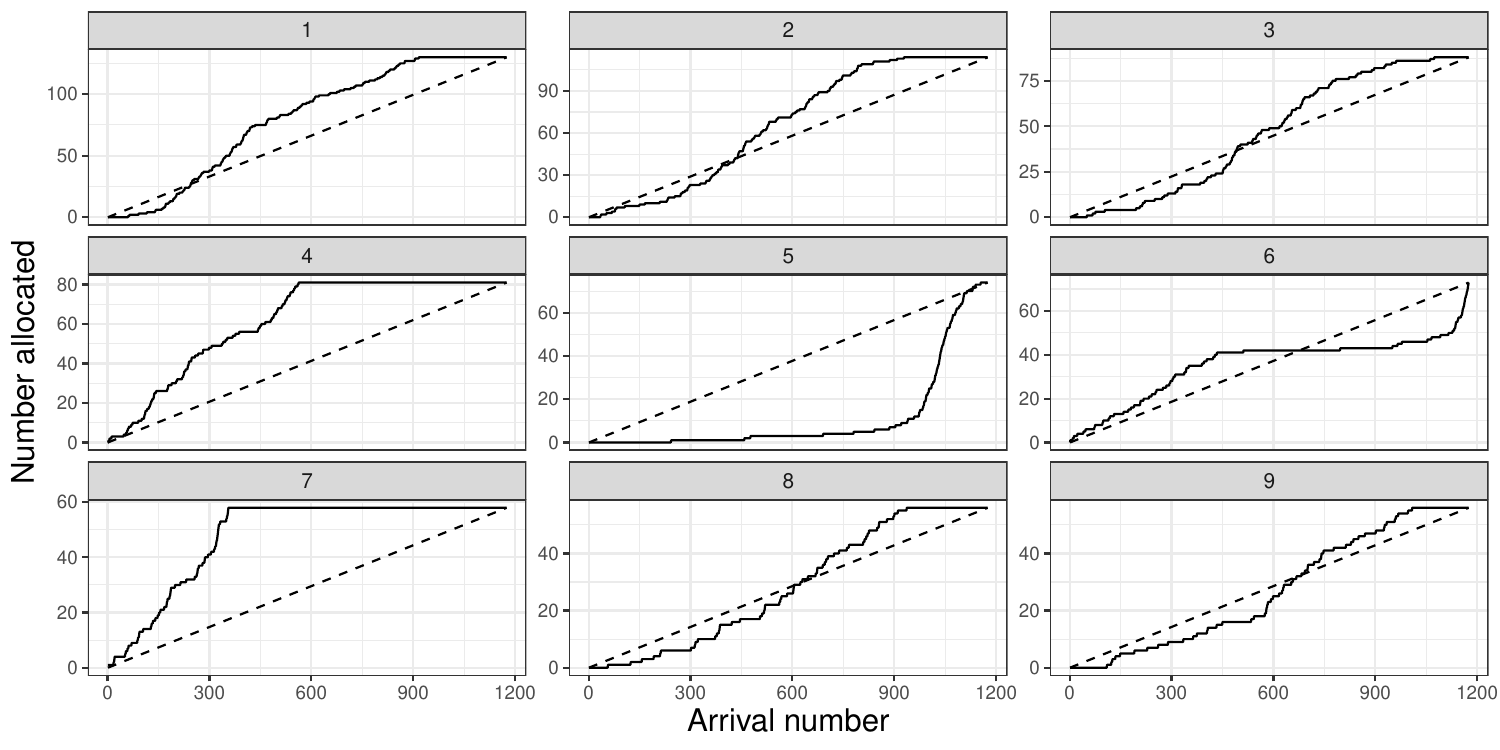}
	\includegraphics[width=.45\textwidth]{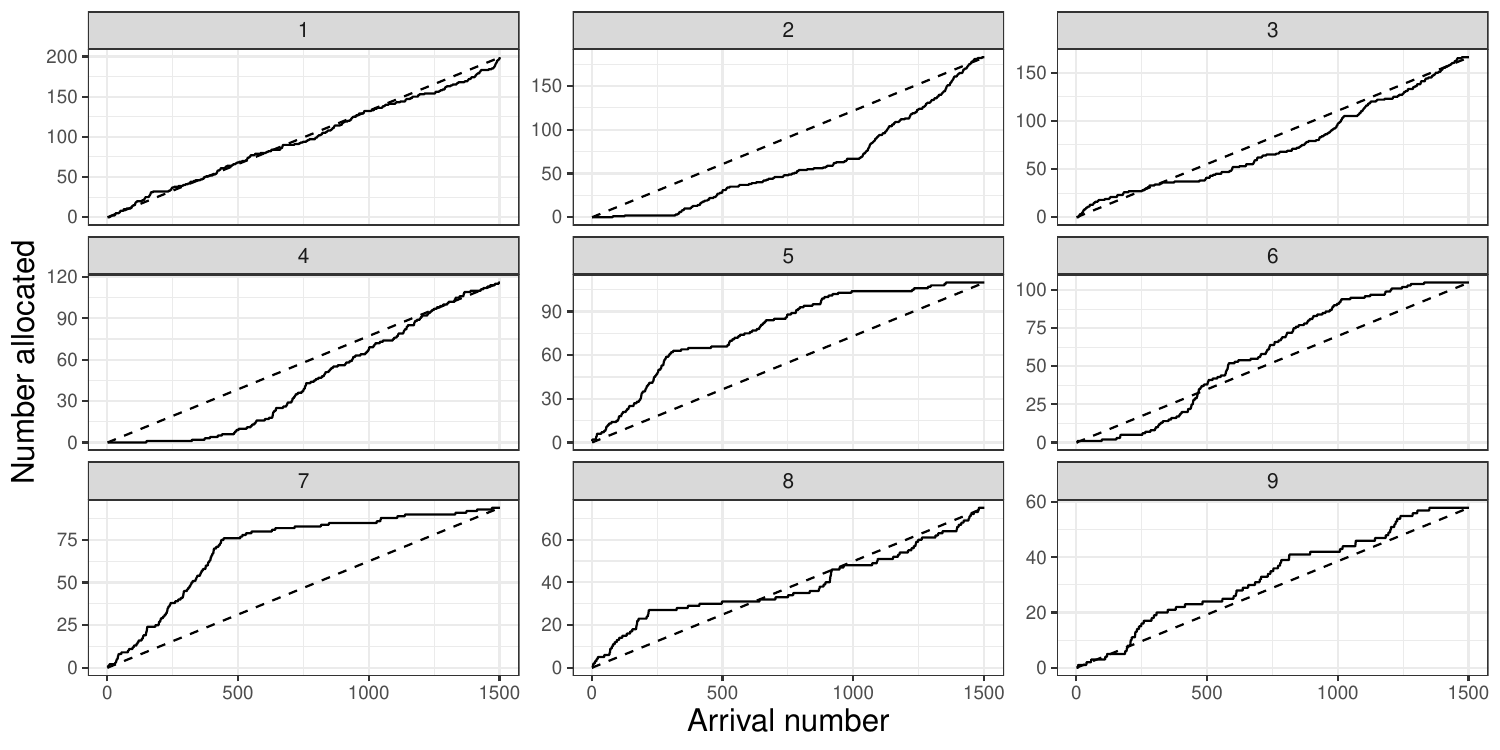}
	\caption{The allocation over time obtained under the potentials method on US data (left) and Swiss data (right). The average queue length on the US data is 7.54, and on the Swiss data is 6.35.}\label{fig:potentials_comparison_imbalance}
\end{figure}

Figure \ref{fig:potentials_comparison} compares the average employment level achieved by the algorithms proposed in this paper to the potentials algorithm of \cite{ahani2021dynamic} on the US data (left) and Swiss data (right). On both datasets, the algorithms result in virtually identical average employment. Figure \ref{fig:potentials_comparison_imbalance} shows the imbalance resulting from the potentials method of \cite{ahani2021dynamic} on the US data (left) and Swiss data (right), which is extremely similar to the imbalance resulting from \ref{alg:modal}. 

In general, we would expect both approaches to perform similarly in terms of solution quality and computational cost, as they can be thought of as primal/dual versions of each other. If $K=1$ and each method samples the same single trajectory, their solution would be identical. When $K>1$, the approaches differ in how results are aggregated across the sample trajectories. In the potentials method, the average of the dual variables across the $K$ trajectories is used to inform the current assignment decision. In \ref{alg:modal}, the most commonly chosen assignment location across the $K$ trajectories is chosen. In terms of run-time, the potentials method requires solving $2K$ LPs for each arrival, each with $M+T$ decision variables. \ref{alg:modal}, on the other hand, requires solving $K$ LPs for every arrival, each with $M\times T$ decision variables. 

\begin{figure}[ht!]
	\centering
	\includegraphics[width=.45\textwidth]{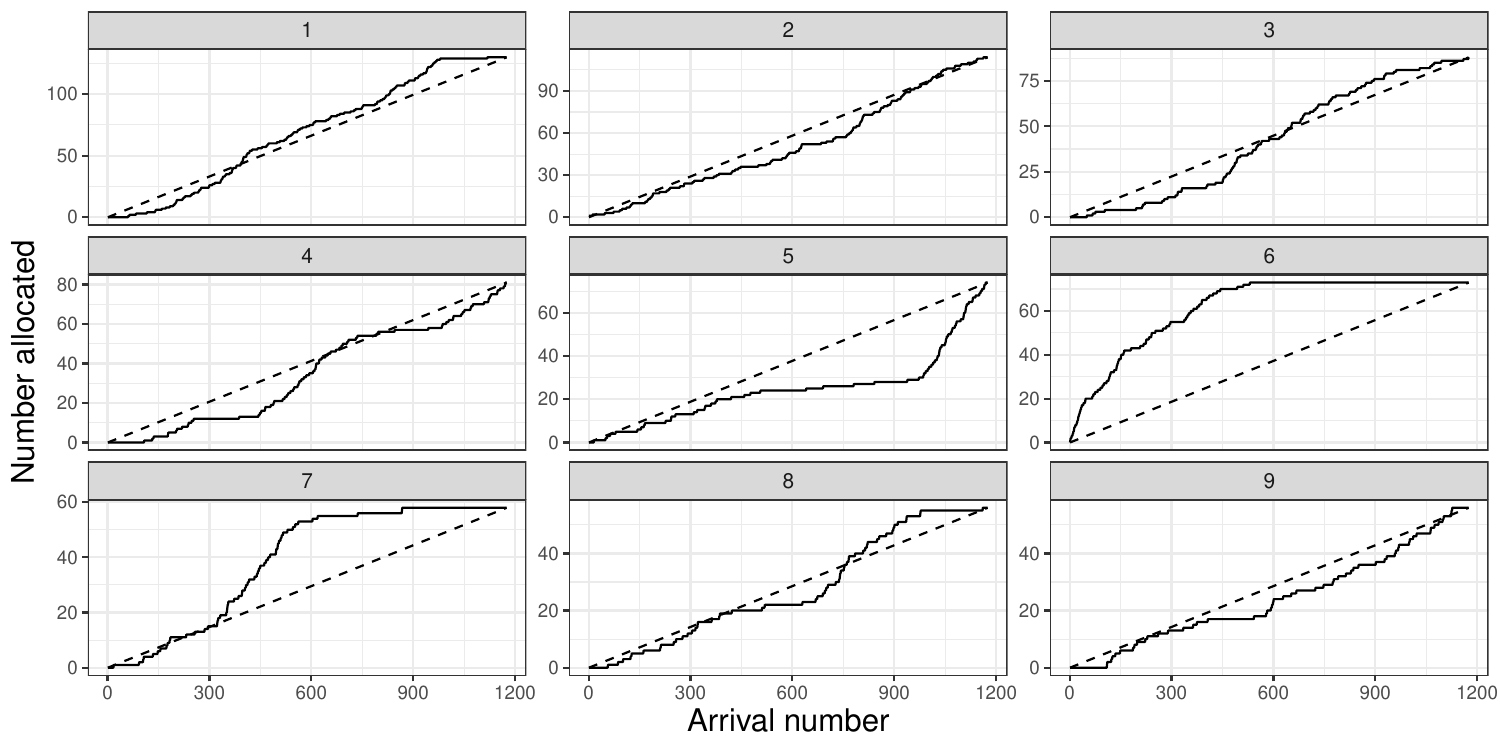}
	\includegraphics[width=.45\textwidth]{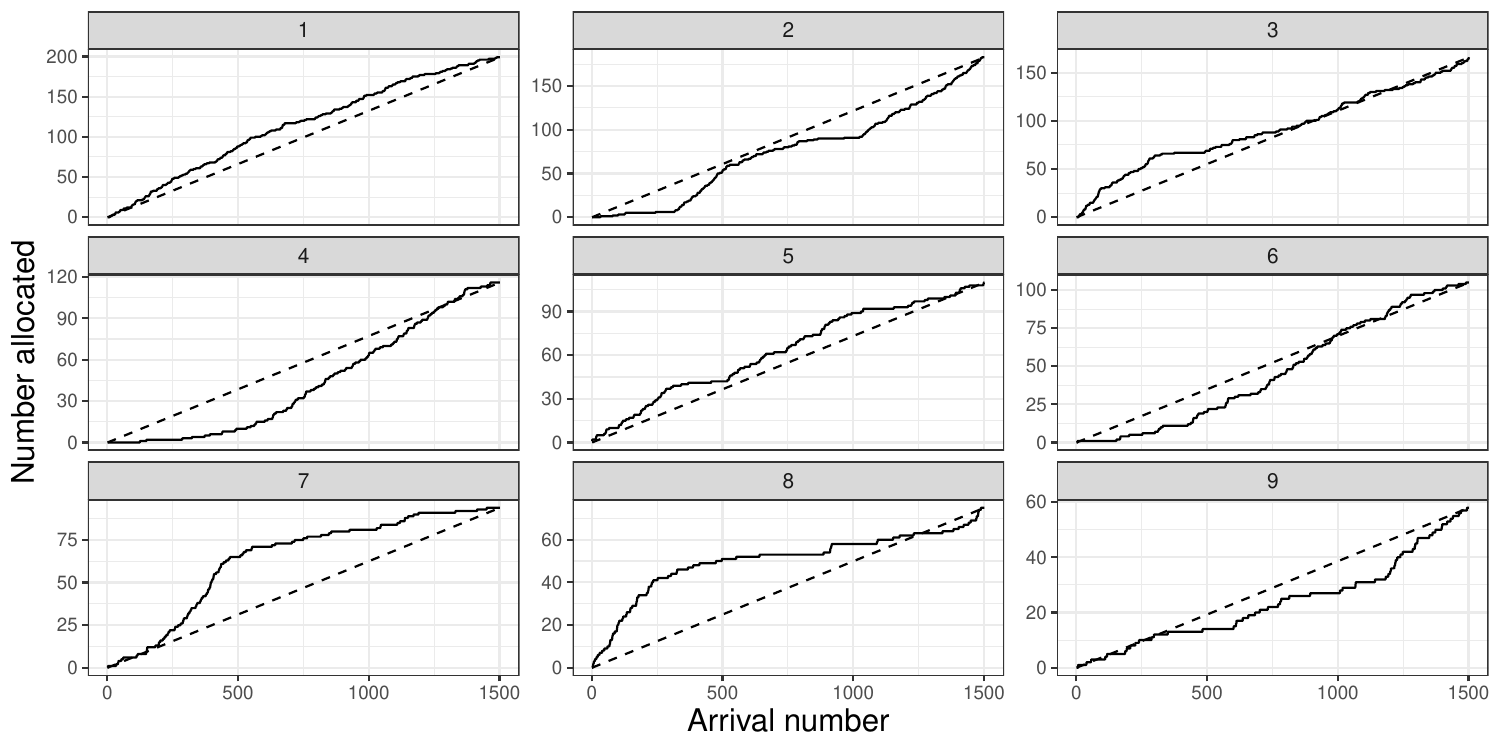}
	\caption{Imbalance of the hindsight-optimal solution for the employment-maximizing assignment problem on US data (left) and Swiss data (right). The average queue length on the US data is 5.30, and on the Swiss data is 5.72.}\label{fig:imbalance_opt}
\end{figure}

Figure \ref{fig:imbalance_opt} shows the allocation to the top nine locations over time resulting from the hindsight-optimal solution. Although the imbalance is not as severe as it is under \ref{alg:modal}, significant imbalance still persists. This further motivates the need to develop an algorithm that explicitly takes balance into account. 

\section{Outcomes for Subgroups}\label{sec:subgroup-outcomes}

\begin{figure}[ht!]
	\centering
	\includegraphics[width=\textwidth]{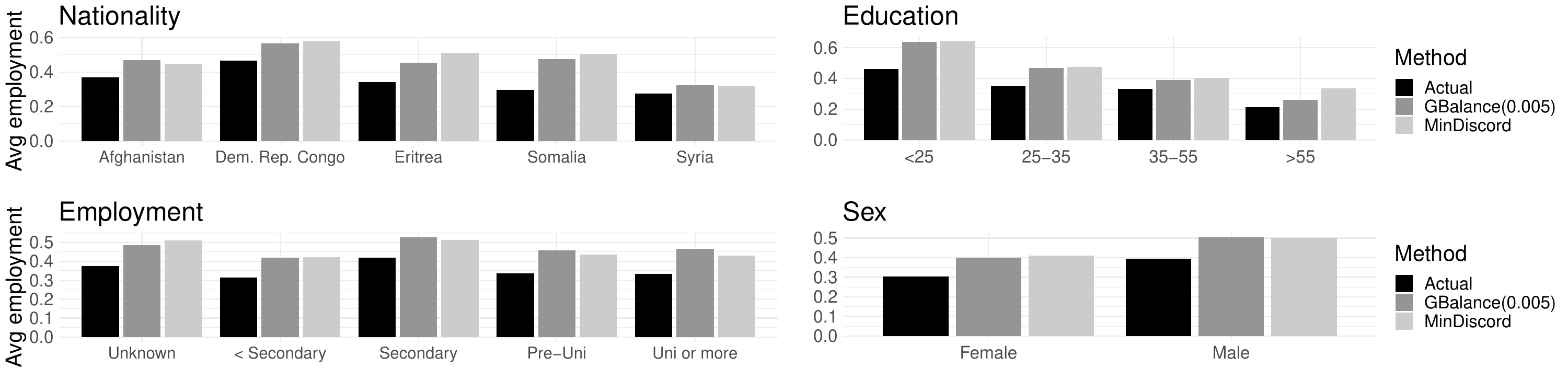}
	\caption{Employment outcomes for subgroups under \ref{alg:modal}, \ref{alg:modal_balance}, and the actual historical assignments. For nationality, only the top 5 most frequent nationalities are shown.}\label{fig:subgroups}
\end{figure}

In this section we report the employment outcomes of the proposed methods for four key subpopulations defined by: sex, age, education, and nationality for the US data. We show the resulting average employment achieved by 1) \ref{alg:modal}, 2) \ref{alg:modal_balance}, and 3) the actual historical assignment. As can be seen in Figure \ref{fig:subgroups}, all subgroups have higher predicted employment rates under the proposed algorithms than under the actual historical assignment. 

\section{Minimum-Risk Formulation} \label{sec:min-risk}

\ref{alg:modal} seeks to minimize the disagreement probability at each time step. An alternative algorithm could instead attempt to directly minimize the expected disagreement cost, which is again given by:
$$d_t = \mathbb{E}_{\mathbf{W}_{t+1}}[\Psi^*_t(\tilde{\mathbf{s}}(t-1))-\left(w_{t \phi(t)} + \Psi_{t+1}^*(\tilde{\mathbf{s}}(t-1)-e_{\phi(t)})\right)  | \mathbf{w}_t]$$
where  $\Psi_{t}^*(\tilde{\mathbf{s}}):=\mathbb{E}_{\mathbf{W}_{t}}[v(\ref{prob:offline}(\mathbf{W}_{t}, \tilde{\mathbf{s}}))]$ denotes the expected sum of employment scores of the optimal assignment for all cases $t$ onward, starting with capacity vector $\tilde{\mathbf{s}}$. 

Algorithms such as the potentials method proposed by \cite{ahani2021dynamic} attempt to minimize the disagreement cost by using the dual variables from the capacity constraints of Problem \ref{prob:offline} as approximations of the following quantity for each location $j$:
$$\Psi^*_t(\tilde{\mathbf{s}}(t-1))-\Psi_{t+1}^*(\tilde{\mathbf{s}}(t-1)-j) $$
This approach (as demonstrated in Figure \ref{fig:potentials_comparison}) appears to perform quite similarly to, if not slightly worse than, \ref{alg:modal} on the data in this study.

The following steps provide a more exact and direct approach to minimizing the disagreement cost in each time step. We refer to this as the ``minimum-risk" approach. \\
For each location $j\in[M]$:
\begin{enumerate}
	\item Sample $K$ arrival trajectories for cases $t+1$ through $T$, denoted by $\mathbf{W}^k_{t+1}$, and compute $d_{tj}^k := \Psi^*_t(\mathbf{W}^k_{t+1},\tilde{\mathbf{s}}(t-1))-\left(w_{t j} + \Psi_{t+1}^*(\mathbf{W}^k_{t+1},\tilde{\mathbf{s}}(t-1)-j) \right) $ for each trajectory and each $j \in [M]$.
	\item Assign case $t$ to $\argminA_j \frac{\sum_{k=1}^Kd_{tj}^k  }{K}$.
\end{enumerate}
Although this method approximates minimizing the disagreement cost exactly, it requires solving $2\cdot K\cdot M$ instances of \ref{prob:offline} in each time step. If run-time is not an issue, or if both the number of locations and horizon length are average in size, this method would be an appropriate choice over \ref{alg:modal}.

\end{document}